\outer\def\give#1. {\medbreak
             \noindent{\bf#1. }}                     
\outer\def\section #1\par{\bigbreak\centerline{\S
     {\bf#1}}\nobreak\smallskip\noindent}
\def\({\left(}
\def\){\right)}

\def\sqr#1#2{{\vcenter{\hrule height.#2pt              
     \hbox{\vrule width.#2pt height#1pt\kern#1pt
     \vrule width.#2pt}
     \hrule height.#2pt}}}
\def\square{\mathchoice\sqr{5.5}4\sqr{5.0}4\sqr{4.8}3\sqr{4.8}3}
\def\qed{\hskip4pt plus1fill\ $\square$\par\medbreak}




\def\cB{{\cal B}}
\def\cC{{\cal C}}
\def\cD{{\cal D}}

\def\cJ{{\cal J}}

\def\cL{{\cal L}}

\def\cN{{\cal N}}

\def\cP{{\cal P}}

\def\cR{{\cal R}}
\def\cS{{\cal S}}
\def\cT{{\cal T}}
\def\cU{{\cal U}}
\def\cV{{\cal V}}
\def\cW{{\cal W}}


\def\CC{{\rm\kern.24em\vrule width.02em height1.4ex depth-.05ex\kern-.26em C}}
\def\RR{{\,\rm{\vrule width.02em height1.55ex depth-.07ex\kern-.3165em R}}}

\def\C{{\bf C}}
\def\cx#1{{\C}^{#1}}     
\def\cp1{{{\bf P}^1}}
\def\R{{\bf R}}

\def\Z{{\bf Z}}
\def\Q{{\bf Q}}
\def\bar{\overline}              




\input epsf.sty
\magnification\magstep1
\overfullrule0pt
\centerline{Polynomial Diffeomorphisms of $\C^2$.}
\centerline{VIII: Quasi-Expansion}
\bigskip
\centerline{Eric Bedford* and John Smillie\footnote*{Research supported in
part by the
NSF.}}
\bigskip
\section 0.  Introduction

This paper continues our investigation of the dynamics of polynomial
diffeomorphisms of $\C^2$ carried out in [BS1-7]. There are several reasons
why the polynomial diffeomorphisms of $\C^2$ form an interesting family of
dynamical systems.  Not the least of these is the fact that  there are
connections with two other areas of dynamics: polynomial maps of $\C$ and
diffeomorphisms of
$\R^2$, which have each received a great deal of attention.
The fact that these three areas are linked makes it interesting to understand
different dynamical notions in these three contexts. One of the
fundamental ideas in dynamical systems is hyperbolicity. One lesson from the
study of the dynamics of maps of $\C$ is that hyperbolicity does not stand
alone
as a dynamical property, rather, it is one of a sequence of interesting
properties which can be defined in terms of recurrence properties of critical
points. These one dimensional properties include the critical finiteness
property,
semi-hyperbolicity, the Collet-Eckmann property and others.  In this paper we
introduce a dynamical property of polynomial diffeomorphisms that generalizes
hyperbolicity  in the way that semi-hyperbolicity generalizes hyperbolicity for
polynomial maps of
$\C$.

In one dimensional complex dynamics generalizations of hyperbolicity are
typically defined in terms of recurrence properties of critical points. Since
we are dealing with diffeomorphisms of $\C^2$ there are no critical points,
and we must use other methods. One way to approach expansion
properties is via a certain canonical metric on unstable tangent spaces of
periodic saddle points that we define.  A mapping is said to the {\it
quasi-expanding} if this metric is uniformly expanded. Although this metric is
canonical it need not be equivalent to the Euclidean metric.  It follows that
quasi-expansion need not correspond to uniform expansion in the usual sense. We
will see in fact that quasi-expansion is strictly weaker than uniform
expansion.
If both $f$ and $f^{-1}$ are quasi-expanding we say that $f$ is
quasi-hyperbolic.  We will show in this paper that quasi-hyperbolic
diffeomorphisms have a great deal of interesting structure.
Using this structure we develop criteria for showing that certain
quasi-expanding diffeomorphisms are uniformly hyperbolic. This criteria for
hyperbolicity (as well as the general structure of quasi-hyperbolic
diffeomorphisms) plays and important role in the study of real diffeomorphisms
of maximal entropy which is carried out in [BS].

We now define the canonical metric on which the definition of
quasi-expansion is
based. We let
$\cS$ denote the set of saddle points of $f$, and we let $J^*$ denote
the closure of $\cS$.  For $p\in\cS$ we let $W^u(p)$ denote the unstable
manifold through $p$, and we let $E^u_p$ denote its
tangent space at $p$.  $W^u(p)$ has the structure of a Riemann
surface immersed in $\C$, and there is a conformal uniformization
$\psi_p:\C\to W^u(p)$ with the property that $\psi_p(0)=p$.  We may
normalize $\psi_p$ by the condition that
$\max_{|\zeta|\le1}G^+\circ\psi_p(\zeta)=1$.  We may define a norm
$\Vert\cdot\Vert^\#$ on $E^u_p$ by the condition that the differential of
$\psi_p$ with respect to the euclidean metric on $\C$ has norm 1 at the
origin.  A mapping is said to the {\it quasi-expanding} if this metric is
expanded by a constant $\kappa>1$ for all $p\in\cS$.

In Section
1 we describe several conditions which are equivalent to quasi-expansion.  One
such condition is that the family of uniformizations $\{\psi_p:p\in\cS\}$ is a
normal family of entire functions. Quasi-expansion is a property of
diffeomorphisms. In Section 2 we consider related properties of individual
orbits.

Let $\Psi$ denote the set of normal limits of $\{\psi_p:p\in\cS\}$, and let
$\cW^u=\{\psi(\C):\psi\in\Psi\}$.  In Section 3 we consider $\psi(\C)$
purely as a variety, that is to say without regard to its parametrization.
For fixed $r>0$ we let $B(p,r)$ denote the ball in $\C^2$ with center $p$ and
radius $r$.  We let $W^u(p,r)$
denote the connected component of
$B(p,r)\cap W^u(p)$ containing $p$.  If $f$ is quasi-expanding, then there
exists
$r>0$ such that for all $p\in\cS$, $W^u_r(p)$ is closed in $B(p,r)$, and
the area of
$W^u(p,r)$ is bounded above.  By Bishop's Theorem and Lemma 2.6,
the correspondence $\cS\ni p\mapsto W^u(p,r)$ extends to a continuous family of
varieties $J^*\ni x\mapsto V^u(x,r)$ such that $V^u(p,r)=W^u(p,r)$ for
$p\in\cS$.  We
prove a ``Bounded Area Distortion Theorem" for proper holomorphic mappings
of planar
domains into $\C^n$ (Theorems 3.1--2).  This is used to prove Theorem 3.4,
which says
that the locally bounded area condition, together with a generalized
transversality
condition, imply quasi-expansion. The bounded area condition also allows us
to prove
that uniform hyperbolicity implies quasi-expansion.

The metric that we define is canonical, but it is not the only canonical
metric that can be defined.  In Section 4, we consider various methods of
defining metrics on unstable tangent spaces $E^u_p$ for $p\in\cS$.  We consider
the equivalence of uniform expansion for various choices of metrics.  In
particular, we define the metrics $\Vert\cdot\Vert^{(L)}$ and
show that they are uniformly expanded by real mappings with maximal
entropy. We show that the uniform expansion of this metric implies
quasi-expansion.

For $x\in J^*$ we let $\Psi_x$ denote the maps $\psi\in\Psi$ with
$\psi(0)=x$.  Such a map has the form
$\psi(\zeta)=x+a_j\zeta^j+O(\zeta^{j+1})$,
and we define the order of $\psi$ to be $j$.  We use the notation $\tau(x)$ for
the maximum order for a function $\psi\in\Psi_x$.  We let $\cJ_j=\{x\in
J^*:\tau(x)=j\}$.  For
$x_0\in\cJ_1$, every function $\psi\in\Psi$ has nonvanishing differential at
the origin.  Thus $\cW^u$ is a lamination in a neighborhood of every point of
$\cJ_1$.  (In Section 6,
$\cJ_1$ will be shown to be a dense, open subset of $J^*$.)  In Section 5,
we show that
$\tau$ describes the local folding of $\cW^u$.  In particular,  $\cW^u$ is
not a
lamination in the neighborhood of $x_1$ if $\tau(x_1)>1$.

In Section 6, we define a metric $\Vert\cdot\Vert^\#_x$ at all points $x\in
J^*$.  This metric (in general not equivalent to the Euclidean metric) is
uniformly expanded if $f$ is quasi-expanding. It follows (Theorem 6.2) that
the largest Lyapunov exponent of a quasi-expanding mapping with respect to
any ergodic invariant measure is strictly positive. In particular it follows
(Corollary 6.3) that all periodic points in $J^*$ are saddle points and
that the
Lyapunov exponents of periodic orbits are uniformly bounded away from 0.

Starting with Section 7 our work applies to mappings for which both $f$ and
$f^{-1}$ are quasi-expanding.  Regularity of the variety $V_x$ at $x$ is shown
for points $x\in\cJ_j$ such that $\alpha(x)\cap\cJ_j\ne\emptyset$.  We also
show that a tangency between $\cW^u$ and $\cW^s$ at a point $x\in J^*$ causes
$\tau(\hat x)>1$ for $\hat x\in\alpha(x)$.

In Section 8 we examine uniform hyperbolicity more carefully.  In
Theorem 8.3 we show that there are geometric properties of $J^\pm$ which
imply hyperbolicity.  In Theorem 8.3, we show that: If
$f,f^{-1}$ are quasi-expanding, and if $f$ is topologically expansive, then
$f$ is uniformly hyperbolic.  Finally, if $f$ and $f^{-1}$ are both
quasi-expanding, we define
the singular set
$\cC$ to be the points of $J^*$ where $\max(\tau^s,\tau^u)>1$.  Let  $\cT$
denote the points of tangency between
$\cW^s$ and
$\cW^u$.  In Theorem 8.10 we show that if $f$ is a mapping for which $\cC$ is
finite and nonempty, then
$\cT\subset W^s(\cC)\cap W^u(\cC)$, and $\bar\cT-\cT=\cC$.

One of our motivations in studying quasi-expansion was to develop the
2-dimensional
analogue of semi-hyperbolicity. In the Appendix, we present the
1-dimensional analogue
of quasi-expansion and show that it is equivalent to semi-hyperbolicity.

\section 1.  Normal Families of Immersions

We say that a holomorphic map $\phi:\C\to\C^2$ is an {\it injective immersion}
(or simply immersion if no confusion will result) if it is injective and an
immersion, which means that $\phi'(\zeta)\ne0$ for all
$\zeta\in\C$.  In this Section we explore the condition that a set of
immersions
has uniform expansion; we show that, in the language of function theory,
this is
equivalent to the set of immersions being a normal family of entire functions.
Let $S\subset J^*$ be a dense, $f$-invariant set.  Suppose that for each
$x\in S$ there is a holomorphic immersion $\psi_x:\C\to\C^2$ such that
$$x\in\psi_x(\C),{\rm\ and\ \ }\psi_x(\C)\subset J^-.\eqno(1.1)$$
In addition, suppose that the family of sets $\{\psi(\C):x\in S\}$
is $f$-invariant, i.e.\
$$f(\psi_x(\C))=\psi_{fx}(\C)\eqno(1.2')$$
and satisfies: for $x_1,x_2\in S$, either $\psi_{x_1}(\C)$ and
$\psi_{x_2}(\C)$ are either disjoint, or they coincide, i.e.\
$$\psi_{x_1}(\C)\cap\psi_{x_2}(\C)\ne\emptyset
\ \ \Rightarrow\ \ \psi_{x_1}(\C)=\psi_{x_2}(\C).\eqno(1.2'')$$
For any holomorphic immersion
$\phi:\C\to\C^2$ with $\phi(\C)=\psi_x(\C)$, there are constants
$a,b\in\C$, $a\ne0$ such that
$\phi(\zeta)=\psi_x(a\zeta+b)$.  We may choose $a$ and $b$ to obtain the
normalization properties:
$$\psi_x(0)=x,\quad \max_{|\zeta|\le1}G^+(\psi_x(\zeta))=1.\eqno(1.3)$$
The first condition in (1.3) may be achieved by a translation of
$\zeta$.  To see that the second normalization is always possible, we note
that since $x\in J^*$, $G^+(x)=G^+(\psi_x(0))=0$.  Thus
$$m_x(r):=\max_{|\zeta|\le r}G^+(\psi_x(\zeta))$$
satisfies $m_x(0)=0$ and is a continuous, monotone increasing function
which is unbounded above.  So after a scaling of $\zeta$, we will have
$m_x(1)=1$.  We note that this normalization defines the parametrization
of $\psi_x$ uniquely, up to replacing $\zeta$ by a rotation
$e^{i\theta}\zeta$,
$\theta\in\R$.

By $\psi_S=\{\psi_x:x\in S\}$ we denote the family of these immersions,
normalized by (1.3).  For $x\in S$ there is a linear mapping
$L_x:\C\to\C$, $L_x(\zeta)=\lambda_x\zeta$, and the family $\{L_x:x\in
S\}$ has the property
$$f\circ \psi_x=\psi_{fx}\circ L_x.\eqno(1.4)$$
Changing the parametrization of $\psi_x$ or $\psi_{fx}$ by a rotation
induces a rotation on
$L_x$.  $f$ induces a mapping $\tilde f$ of $\psi_S$ to itself, given by
$$\tilde f(\psi_x)=\psi_{fx}=f\circ \psi_x\circ L_x^{-1}.\eqno(1.5)$$
For $n>0$ set
$$\lambda(x,n)=\lambda_x\lambda_{fx}\dots\lambda_{f^{n-1}x}\eqno(1.6)$$
for $n>0$, then
$$\tilde
f^n(\psi_x)(\zeta)=f^n\circ\psi_x(\lambda(x,n)^{-1}\zeta).\eqno(1.7)$$

By
the identity
$G^+\circ f=d\cdot G^+$ and the transformation formula $(1.2')$, we have
$$d\cdot m_x(r)=m_{fx}(|\lambda_x|r).\eqno(1.8)$$
Setting $r=1$, we have
$d=m_{fx}(|\lambda_x|)>1$, which gives $|\lambda_x|>1$ for all $x\in S$.

For $x\in S$, we let $E^u_x$ denote the subspace of the tangent space of
$T_x\C^2$ given by the $\C$-linear span of $\psi_x'(0)$.  If $v\in
E^u_x$, then $v$ is a scalar multiple of
$\psi_x'(0)$, so we we define the norm
$$\|v\|_x^\#:=|v/\psi_x'(0)|.\eqno(1.9)$$
It
follows that the norm of
$Df_x$, measured with respect to this family of norms, is given by
$$\|Df_x|_{E^u_x}\|^\# =\max_{v\in E^u_x-\{0\}}{\|Df_xv\|_{fx}^\#\over
\|v\|_x^\#}=|\lambda_x|.$$
Similarly, $\|Df_x^n|_{E^u_x}\|^\#=|\lambda(x,n)|$.

For $1\le r<\infty$ we define
$$M(r)=\sup_{x\in S}m_x(r).$$
For each $x\in S$,  $m_x(r)$ is a convex, increasing function
of $\log r$.  It follows that $M(r)$ also has these properties on the open
interval where it is finite.  In particular, $M(r)$ is continuous from the
right at $r=1$, $M(1)=1$, and $M(r)>1$ if $r>1$.

\proclaim Lemma 1.1.  $\psi_S$ is a normal family if and only if
$M(r)<\infty$ for all $r<\infty$.

\give Proof.  We set $V=\{|x|,|y|\le R\}$ and
$V^+=\{|y|\ge|x|,|y|\ge R\}$.  It is known that for $R$ sufficiently large,
$W^u(x)\subset V\cup V^+$, $S\subset V$, and
$G^+|_{V\cup V^+}$ is a proper exhaustion.  Since $\psi_x(0)=x\in S$, it
follows that no sequence in $\psi_S$ can diverge to infinity uniformly on
compacts.  Thus normality is equivalent to local boundedness at every
point.  For fixed
$\zeta\in\cx{}$ the sequence $\{\psi_{x_j}(\zeta)\}$ is bounded if
and only if $\{G^+(\psi_{x_j}(\zeta))\}$ is bounded.  Since $M(r)$ is
increasing in $r$, it follows that if $M(r)<\infty$, then $\Psi$ is a
normal family on
$\{|\zeta|<r\}$.

Conversely, if $\psi_S$ is a normal family, then
$\{G^+(\psi(\zeta)):\psi\in\Psi,|\zeta|\le r\}$ is bounded.  Thus
$M(r)<\infty$.  \qed

The following result shows that the normal family condition is equivalent to
a number of ``uniform conditions''.

\proclaim Theorem 1.2.  The following are equivalent:
\item{(1)} $\psi_S$ is a normal family.
\item{(2)} $M(r_0)<\infty$ for some $1<r_0<\infty$.
\item{(3)} For all
$r_1<r_2<\infty$ there is a constant $k<\infty$ such that $m_x(r_2)/m_x(r_1)\le
k$ for all $x\in S$.
\item{(4)} There exists $\kappa>1$ such that for all $x\in S$,
$|\lambda_x|\ge\kappa$.
\item{(5)} ${\rm There\ exist\ }C,\beta<\infty{\rm \ such\ that\ } m_x(r)\le
Cr^\beta{\rm\ for\ all\ }x\in S{\rm\ and\ }r\ge1. \ \ \ \ \ \ \ \ \ \ \ \
\hfill(1.10)$

\give Proof.  $(1)\Rightarrow(2)$ is a consequence of Lemma 1.1.

$(2)\Rightarrow(4)$  If $M(r_0)<\infty$ for some $1<r_0<\infty$ then $\log
M(r)$ is a convex increasing function of $\rho=\log r$ on the interval
$\rho\in(0,\log r_0)$.  It follows that $M(r)$ is continuous at
$r=1$.  Thus $\kappa:=\inf\{t\ge1:M(t)\ge d\}>1$.  Now for any
$x\in S$ we have
$m_x(|\lambda_x|)=d\cdot m_{fx}(1)=d$.  It follows, then, that
$M(|\lambda_x|)\ge d$, and so
$|\lambda_x|\ge \kappa>1$.

$(4)\Rightarrow(5)$  For $x\in S$, let $x_j=f^jx$.  Then by the
transformation formula (1.8)
$$m_{x}(\kappa^p)\le
m_{x}(|\lambda_{x_0}\lambda_{x_{-1}}\cdots\lambda_{x_{-p+1}}|)=d^p
m_{x_{-p+1}}(1)=d^p.$$ For any $1\le r<\infty$ we choose $p$ such that
$\kappa^{p-1}\le r<\kappa^p$.  If we choose $\beta=\log
d/\log\kappa$, then $\kappa^\beta=d$, and
$$m_x(r)\le m_x(\kappa^p) \le (\kappa^p)^\beta=\kappa^\beta
r^\beta.$$
Thus (1.10) holds with
$C=\kappa^\beta$.

$(5)\Rightarrow(1)$  Condition (1.10) implies that $M(r)\le Cr^\beta$, and
thus $\psi_S$ is a normal family by Lemma 1.1.

$(1)\Leftrightarrow(3)$  Let $\tilde\psi_S$ denote the set of scaled
functions $\tilde\psi(\zeta)=\psi(r_1\zeta)$ for all $\psi\in\psi_S$.  By
the equivalence (1)$\Leftrightarrow$(2) and Lemma 1.1, we have that
$\tilde\psi_S$ is a normal family if and only if
$$\tilde M(r):=\sup_{\tilde\psi\in\tilde\psi_S}{\tilde m_x(r)\over\tilde
m_x(1)}=\sup_{\psi\in\psi_S}{m_x(rr_1)\over m_x(r_1)}<\infty.$$
Finally, it is evident that $\tilde\psi_S$ is a normal family if and only
if $\tilde\psi_S$ is normal.  Thus (1) is equivalent to (3). \qed

We say that $f$ is {\it quasi-expanding} if the equivalent conditions in
Theorem 1.2 hold.  While these conditions are stated in terms of the
family $\psi_S$, we will see in \S3 that they are independent of the
choice of the particular family $\psi_S$.  We say that $f$ is {\it
quasi-contracting} if
$f^{-1}$ is quasi-expanding.

\proclaim Proposition 1.3.  For $n\ge1$, $f$ is quasi-expanding if and
only if $f^n$ is quasi-expanding.

\give Proof.  Let $\psi_S$ be a family satisfying (1.1--3) for $f$.  For
$n\ge1$,
$ J^*$ and $K^+$ are the same for $f^n$.  It follows that $\psi_S$ also
satisfies (1.1--3) for $f^n$.  If $f$ is quasi-expanding, then $\psi_S$ is
a normal family; thus $f^n$, too, is quasi-expanding.

Now suppose that $\psi_S$ satisfies (1.1--3) for $f^n$.  It follows
that $\tilde S:=S\cup fS\cup\dots\cup f^{n-1}S$ is $f$-invariant.  Let
$\psi^{(j)}_S$ denote the set of mappings $\{f^j\circ\psi\circ
L_\lambda:\psi\in\psi_S\}$, where $L_\lambda(\zeta)=\lambda(\psi,j)\zeta$
is chosen so that $f^j\circ\psi\circ L_\lambda$ satisfies the normalization
(1.3).  Let $\tilde\psi_S:=\psi^{(0)}_S\cup\dots\cup\psi_S^{(n-1)}$, so
that $\tilde\psi_S$ satisfies (1.1--3) for $f$.  Define $M^{(j)}(r)
=\sup_{\psi\in\psi^{(j)}_S}\sup_{|\zeta|\le r}G^+\circ\psi(\zeta)$.  If
$f^n$ is
quasi-expanding, then $\psi^{(0)}_S$ is a normal family.  By Lemma 1.1 this
means that
$M^{(0)}(r)<\infty$ for $r<\infty$.  As in the line following (1.8) we have
$\lambda(\psi,j)|\ge1$.  It follows that
$$\eqalign{ M^{(j)}(r)&=\sup_{\psi\in\psi^{(j)}_S}\sup_{|\zeta|\le
r}G^+\circ\psi(\zeta)
\cr &= \sup_{\psi\in\psi^{(0)}_S}\sup_{|\zeta|\le r} d\cdot
G^+\circ\psi(|\lambda(\psi,j)|^{-1}\zeta) \le d\cdot M^{(0)}(r).\cr}$$
Thus $M^{(j)}(r)<\infty$ for all $r<\infty$.  It follows from Lemma 1.1
that each
$\psi^j_S$ is a normal family.  Thus $\tilde\psi_S$ is normal, and $f$ is
quasi-expanding.\qed

\proclaim Proposition 1.4.  If $f$ is quasi-expanding, then for $x\in
S$, $\psi_x(\C)\subset W^u(x)$, i.e.\ if $y_1,y_2\in\psi_x(\C)$, then
$\lim_{n\to+\infty}dist(f^{-n}y_1,f^{-n}y_2)=0$.

\give Proof.  For $j=1,2$ there exist $\zeta_j\in\C$ such that
$\psi_x(\zeta_j)=y_j$.  By (1.7) $f^{-n}y_j=f^{-n}\psi_x(\zeta_j)
=\psi_{x_{-n}}(\lambda(x,-n)^{-1}\zeta_j)$.  Now $\{\psi_{x_{-n}}:n\ge0\}$
is a normal family, so the set of derivatives
$\{|D\psi_{x_{-n}}(\zeta)|:|\zeta|\le1,n\ge0\}$ is uniformly bounded by
$M<\infty$.  Thus
$$\eqalign{
dist(f^{-n}y_1,f^{-n}y_2)=
|&\psi_{x_{-n}}(\lambda(x,-n)^{-1}\zeta_1)
-\psi_{x_{-n}}(\lambda(x,-n)^{-1}\zeta_2)|\cr
&\le|\lambda(x,-n)|^{-1}M|\zeta_1-\zeta_2|,\cr}$$
which tends to zero, since $\lambda(x,-n)\to\infty$ by (4) of Theorem
1.2. \qed

We give two examples to show that families $\psi_S$ satisfying (1.1--3)
exist for any map
$f$.  Let  $p$ be a saddle point, i.e. a periodic point of saddle type.
The stable and unstable manifolds $W^s(p)$ and $W^u(p)$ through $p$ are
conformally equivalent to $\C$.  Let $\phi:\C\to W^u(p)$ denote a
uniformization of the unstable manifold.  It is evident that
$p\in\phi(\C)=W^u(p)$, and by the argument of [BS1, Proposition 5.1] we
have $W^u(p)\subset J^-$.

\give Example 1.  Let $p$ and $q$ be saddle points, and set
$S=W^s(q)\cap W^u(p)$.  By [BLS], $S$ is a dense subset of $ J^*$.  Let
$\phi$ denote the uniformization of $W^u(p)$ as above.  For $x\in
S\subset W^u(p)$, we set $\beta_x:=\phi^{-1}(x)$.  Now we may choose
$\alpha_x\ne0$ such that $\psi_x(\zeta):=\phi(\alpha_x(\zeta+\beta_x))$
satisfies (1.1--3).

\give Example 2.  Let $S$ denote the set of saddle (periodic) points.  By
[BLS], $S$ is dense in $ J^*$.  For $p\in S$ the unstable manifold $W^u(p)$
may be normalized to satisfy the conditions (1.1--3) above.

If $p$ is
periodic of period $n=n_p$, then the multiplier $Df^n_p|_{E^u_p}$ is given
by $\lambda(p,n)$.  Then we have
$$d^{n}m_p(r)=m_p(|\lambda(p,n)|r)$$
Thus we conclude that for $p\in S$
$$m_p(r)\le C_pr^\beta\eqno(1.11)$$
holds with $\beta={\log d/
({1\over n}\log|\lambda(p,n)|)}$.
This condition (1.11) allows both $C$ and $\beta$ to vary with $p$ and
differs from $(1.10)$ in this respect.

A variant of (1.10) is
$${\rm There\ exist\ }C,\beta<\infty{\rm\ such\ that\ \ \
}\inf_{y\in S}m_y(r)\ge{r^\beta\over C}{\rm\ \ \ for\ all\
}0<r<1.\eqno(1.12)$$
\proclaim Proposition 1.5.  If $f$ is quasi-expanding, then (1.12) holds.

\give Proof.  If $f$ is quasi-expanding, then there exists $\kappa>1$ such
that $|\lambda_\psi|\ge\kappa$ for all $\psi\in\Psi$.  Given $r$, choose $n$
such that $\kappa^{-n}\le r<\kappa^{-n+1}$.  By the normalization condition,
$m_\psi(|\lambda(\psi,-n)|)=d^{-n}$.  Thus by the choice of $n$,
$$m_\psi(r)\ge m_\psi(\kappa^{-n})\ge m_\psi(|\lambda(\psi,-n)|)=d^{-n} \ge
d^{{\log
r\over\log\kappa}-1}.$$
Thus (1.12) holds with $C=d$ and $\beta=\log d/\log\kappa$. \qed

\proclaim Corollary 1.6.  If $f$ is quasi-expanding, there are $C>0,
r_0>0$, and
$m<\infty$ such that $\max_{B(x,r)}G^+\ge Cr^m$ for $x\in J^*$ and $0<r<r_0$.

As an alternative to (1.3) we may consider the normalization
$$\psi_x(0)=x,\ \ \ \max_{|\zeta|\le1}G^+\psi_x(\zeta)=t\eqno(1.13)$$
for fixed $0<t<\infty$.  In this case we have a family
$\tilde\psi_S$, normalized according to (1.13).  Further, we have
functions $\tilde m_x(r)$, $\tilde M(r)$, and multipliers
$\tilde\lambda_x$.  We may repeat the proof of Theorem 1.2 and obtain the
equivalence between conditions (1), (2), (3), (5), and the condition (4):
{\sl There exists $\tilde\kappa>1$ such that $|\tilde\lambda_x|>1$ for all
$x\in
S$.}

\give Proposition 1.7.  The family $\tilde\psi_S$, normalized according to
(1.13) is normal if and only if $\psi_S$ is normal.

\give Proof.  By the preceding remarks, it suffices to show that if
$\tilde\psi_S$ is normal, then $\psi_S$ is normal.  If $t\le1$, then
$\tilde M(r)\le M(r)$, and so $\tilde\psi_S$ is normal by Lemma 1.1.  Now
suppose $t>1$.  For each $x\in S$ there is a number $\tau=\tau_x\in\C$,
$|\tau|>1$ such that $\tilde\psi_x(\zeta)=\psi_x(\tau\zeta)$.  Choose $k$
such that $d^k\ge t$.  Then
$$m_x(\kappa^k)\ge m_x(|\lambda(f^{-k}x,k)|) = d^{-k} m_x(1) = d^k\ge t.$$
It follows that $|\tau|\le\kappa^k$, and this upper bound gives the
normality of $\tilde\psi_S$.\qed

\section 2.  Expansion Along Individual Orbits and Unstable Germs

While our primary focus is the dynamics of quasi-expanding
diffeomorphisms, some
of the results in the sequel are local results and depend only on information
about the behavior of a particular orbit.  In this section we explore various
orbitwise notions of expansion and regularity. This section may be omitted
on a first reading of this paper.  We define
$$\eqalign{ &M_x(r)=\limsup_{S\ni y\to x} m_y(r), \ \ \
\hat\lambda_x=\liminf_{S\ni
y\to x}|\lambda_y|\cr
& r_x=\inf\{r:M_x(r)>0\}, \ \ \ R_x=\inf\{r:M_x(r)=\infty\},\cr}$$
where we admit $+\infty$ as a possible value.
It follows that $\hat\lambda_x\ge1$, $M_x(0)=0$, $M_x(1)=1$, and $R_x\ge1$.
Further,
$M_x$ is a convex, increasing function of
$\log r$ for $r$ in the interval $(0,R_x)$, and it is evident that
$x\mapsto M_x$ and
$x\mapsto r_x$ are upper semicontinuous; and
$x\mapsto\hat\lambda_x$ and $x\mapsto R_x$ are lower semicontinuous.

For $n\ge0$, we define
$$\hat\lambda(x,n)=\hat\lambda_x\hat\lambda_{fx}
\dots\hat\lambda_{f^{n-1}x},{\rm\
\ and\ \ \ }
\hat\lambda(x,-n)=\hat\lambda^{-1}_{f^{-n}x}\dots\hat\lambda^{-1}_{f^{-1}x}
=\hat\lambda(f^{-n}x,n)^{-1}.$$

\proclaim Lemma 2.1.  For $x\in J^*$ and $n\ge0$ we have
$$d^nM_{x}(r)\ge M_{f^nx}(\hat\lambda(x,n)r){\rm\ \ and\ \ }d^{-n}M_{x}(r)\le
M_{f^{-n}x}(\hat\lambda(x,-n)r).\eqno(2.1)$$

\give Proof.  By (1.8) we have
$$m_{f^np}(|\lambda(p,n)|r)=d^nm_p(r) {\rm\ \ and\ \ }
m_{f^{-n}x}(r)=d^{-n}m_x(|\lambda(x,-n)|r) $$
  for $p\in S$ and $n\ge0$.
Let us fix $r$ and choose a sequence of points $p_j\to x$ such that
$m_{f^np_j}(\hat\lambda(x,n)r)\to M_{f^nx}(\hat\lambda(x,n)r)$.  Since the
$m_{p_j}$
are convex in $\log r$, then are equicontinuous, so by the lower
semicontinuity of
$x\mapsto\hat\lambda(x,n)$, it follows that
$M_{f^nx}(\hat\lambda(x,n)r)\le\limsup_{p_j\to
x}m_{p_j}(\lambda(p_j,n)r)=\limsup_{p_j\to x}d^nm_{p_j}(r)\le d^nM_{x}(r)$.
The proof
for the other inequality follows by a similar argument, with the only
difference being
that $x\mapsto\hat\lambda(x,-n)$ is upper semicontinuous.
\qed

\proclaim Proposition 2.2.  For $x\in J^*$ and $n\ge0$ we have
$\hat\lambda(x,n)R_x\le
R_{f^nx}$ and $\hat\lambda(x,-n)r_{f^{-n}x}\ge r_x$.  In particular,
$$r_x\le\hat\lambda(x,n)^{-1}\le
1\le\hat\lambda(x,-n)^{-1}\le R_x.\eqno(2.2)$$

\give Proof.  If $r<R_x$, then by Lemma 2.1, we have
$M_{f^nx}(\hat\lambda(x,n)r)<\infty$.  Thus $\hat\lambda(x,n)r\le R_{f^nx}$.
The other inequality is similar. \qed

\proclaim Theorem 2.3.  The following are equivalent:
\item{(1)} $f$ is quasi-expanding
\item{(2)} $\hat\lambda_x>1$ for all $x\in J^*$.
\item{(3)} $\inf_{x\in J^*}\hat\lambda_x>1$.
\item{(4)} $R_x>1$ for all $x\in J^*$.
\item{(5)} $R_x=\infty$ for all $x\in J^*$.
\item{(6)} $\lim_{n\to-\infty}\hat\lambda(x,n)=0$ for all $x\in J^*$.

\give Proof. $(1)\Rightarrow(2)$.  If $f$ is quasi-expanding then
$|\lambda_p|\ge\kappa>1$ for all $p\in S$.  Thus $\hat\lambda_x\ge\inf_{p\in
S}|\lambda_p|\ge\kappa>1$.  $(2)\Rightarrow(3)$.  This follows because
$x\mapsto\hat\lambda_x$ is lower semicontinuous.  $(3)\Rightarrow(1)$.  By the
definition of $\hat\lambda_x$ and the compactness of $J^*$, $\inf_{p\in
S}|\lambda_p|=\inf_{x\in J^*}\hat\lambda_x$.  If
$\kappa:=\inf\hat\lambda_x>1$, then
$|\lambda_p|\ge\kappa$, so $f$ is quasi-expanding.

$(1)\Rightarrow(5)$.  This is condition (2) of Theorem 1.2.

$(5)\Rightarrow(4)$.
This is trivial.

$(4)\Rightarrow(1)$.  Since $x\mapsto R_x$ is lower
semicontinuous,
it follows that $R:=\inf_{x\in J^*}R_x>1$.  Choose $1<R'<R$.  By the upper
semicontinuity of $x\mapsto M_x$, it follows that $\sup_{x\in
J^*}M_x(R')<\infty$.
Thus $f$ is quasi-expanding by (2) of Theorem 1.2.

$(1)\Rightarrow(6)$.  If $f$ is quasi-expanding, then
$\hat\lambda(x,n)\le\kappa^{n}$,
so (6) holds.

$(6)\Rightarrow(5)$. If (6) holds, then $R_x=\infty$ by
Proposition 2.2,
so (5) holds.
\qed

We say that $f$ has {\it forward expansion at $x$} if
$\lim_{n\to+\infty}\hat\lambda(x,n)=\infty$, and we say  that $f$ has {\it
backward contraction at
$x$} if $\lim_{n\to-\infty}\hat\lambda(x,n)=0$. By Proposition 2.2, if
$f$ has forward expansion at $x$, then $r_x=0$; and if $f$ has backward
contraction at $x$, then $R_x=\infty$.

For $x\in J^*$ and $R<R_x$ there is a neighborhood $\cN$ of $x$ in $J^*$ such
that if $y\in \cN$, then $m_y(R)\le M_x(R)+1<\infty$.  Thus the restrictions
$\{G^+\circ\psi_y|\{|\zeta|<R\}:y\in\cN\}$ are uniformly bounded.  Since
$\psi_y(\C)\subset J^-$, it follows that the restrictions of
$\{\psi_y:y\in\cN\}$ to
$\{|\zeta|<R\}$ are uniformly bounded and are thus a normal family.
We let $\Psi_x$ denote the set of analytic mappings
$\psi:\{|\zeta|<R_x\}\to\C^2$
which are obtained as normal limits $\lim_{y_j\to
x}\psi_{y_j}|_{\{|\zeta|<R_x\}}$
for sequences $y_j\to x$.  We set
$\Psi=\bigcup_{x\in J^*}\Psi_x$.  In general it may happen that an element
$\psi\in\Psi_x$ may be analytically extended to a domain strictly larger than
$\{|\zeta|<R_x\}$.  The size of the domain $\{|\zeta|<R_x\}$ assures
that $\Psi$ is a normal family.

Let us define a condition at a point $x$:
$$\Psi_x \hbox{\rm\ contains\ a\ non-constant\ mapping.}\eqno(\dag)$$
Suppose $(\dag)$ holds, and choose a non-constant
$\psi\in\Psi_x$.  We say that $\psi$ is a maximal element of $\Psi_x$ if
whenever
$\psi(\alpha\zeta)$ also belongs to $\Psi_x$ for some constant
$\alpha\in\C$, we have
$|\alpha|\le1$.  By the compactness of $\Psi_x$, each $\psi\in\Psi_x$ is either
maximal or has the form $\psi(\zeta)=\hat\psi(\alpha\zeta)$ for some maximal
$\hat\psi$ and $|\alpha|\le1$.  Passing to convergent subsequences in (1.7)
we see that
if
$\psi\in\Psi_x$ is maximal, then there are a unique (modulo rotation of
variable)
linear transformation
$L(\zeta)=\lambda_\psi\zeta$ and a unique maximal $\psi_1\in\Psi_{fx}$ such
that
$f\circ\psi\circ L^{-1}_\psi=\psi_1$.  This allows us to define
$$\tilde f:\Psi_x\to\Psi_{fx}, \ \ \ \tilde
f(\psi)(\zeta)=f(\psi(\lambda_\psi^{-1}\zeta)).\eqno(2.3)$$
If $\psi$ is not maximal, and if $\psi(\zeta)=\hat\psi(\alpha\zeta)$ is as
above, then
we set $\tilde f(\psi)(\zeta):=\tilde f(\hat\psi)(\alpha\zeta)$.

We use the notation $\psi_j:=\tilde f^j(\psi)$ and
$$\lambda(\psi,n)=
\lambda_{\psi_0}\lambda_{\psi_1}\dots\lambda_{\psi_{n-1}}$$ so
$$\tilde f^n(\psi)(\zeta)=f^n\circ \psi(\lambda(\psi,n)^{-1}\zeta).$$
Since each $\psi\in\Psi$ is a limit of elements of $\psi_S$, and
$\hat\lambda_x$ is a
lim-inf, we have the following.

\proclaim Corollary 2.4.  If (\dag) holds at $x\in J^*$, then for $n\ge0$,
$$\hat\lambda(x,n)=\inf_{\psi\in\Psi_x}|\lambda(\psi,n)|=
\min_{\psi\in\Psi_x}|\lambda(\psi,n)|,$$
where the infimum and minimum are taken over all nonconstant elements of
$\Psi_x$.  In
particular, if $x$ is a point of forward expansion, then
for all $\psi\in\Psi$, $|\lambda(\psi,n)|\to\infty$ and $n\to+\infty$; and
if $x$ is a
point of backward contraction, then $|\lambda(\psi,n)|\to0$ as $n\to-\infty$.

\proclaim Corollary 2.5.  Suppose that (\dag) holds at each $x\in J^*$.  If
$f$ is not
quasi-expanding, then there exists a nonconstant $\psi\in\Psi$ such that
$|\lambda(\psi,n)|=1$ for all $n\le0$.

\give Proof.  If $f$ is not quasi-expanding, then by Theorem 2.3,
$R_x=1$ for some $x\in J^*$.   By Lemma 2.2, we have $\hat\lambda(x,n)=1$
for all
$n\le0$.  By Corollary 2.4, there is a nonconstant $\psi^{\langle
n\rangle}\in\Psi_x$
such that $|\lambda(\psi^{\langle n\rangle},n)|=1$.  By the compactness of
$\Psi_x$, we
may choose a subsequence such that $\psi^{\langle
n_j\rangle}\to\psi\in\Psi_x$, and
$\psi$ has the desired property.
\qed

Let $\psi$ denote the germ at $\zeta=0$ of a nonconstant holomorphic map from a
neighborhood of the origin in $\C$ to $\C^2$.  Setting $x=\psi(0)$, it
follows that
$\{|\zeta|<r,\psi(\zeta)=x\}=\{0\}$ for $r>0$ sufficiently small.  Let
$B(x,\epsilon)$
denote the Euclidean ball in $\C^2$ with center $x$ and radius $\epsilon$, and
let $V(\psi,\epsilon)$ denote the connected component of
$B(x,\epsilon)\cap\psi(|\zeta|<r)$ containing $x$.  If
$\epsilon<\min_{|\zeta|=r}|\psi(\zeta)|$, then
$V(\psi,\epsilon)$ is an analytic subvariety of $B(x,\epsilon)$.

\proclaim Lemma 2.6.  If (\dag) holds, the nonconstant elements of
$\Psi_x$ define a unique germ of a complex analytic variety at $x$.

\give Proof.  Let $\psi_1,\psi_2\in\Psi_x$ be given, and let $V_1,V_2$
be the corresponding germs of varieties, defined in some ball
$B(x,\epsilon)$.  If $\psi_1^j,\psi_2^j$ are sequences from $\psi_S$
which converge to $\psi_1,\psi_2$, respectively, then for $j$
sufficiently large, $\psi_1^j,\psi_2^j$ define subvarieties
$V_1^j,V_2^j$ (respectively) of $B(x,\epsilon)$.  If $V_1$ and $V_2$
define distinct germs of varieties at $x$, then $V_1$ and $V_2$ have a
0-dimensional intersection in $B(x,\epsilon)$.  Thus $V_1^j$ and $V_2^j$
also have 0-dimensional intersection in $B(x,\epsilon)$, which
contradicts $(1.2'')$.
\qed
We will sometimes use the notation $V(x,\epsilon)$ for $V(\psi,\epsilon)$;
and we will let $V_x$ denote the corresponding germ at $x$, which is
independent
of $\psi$ by Lemma 2.6.

We may define
$$\tilde V_x=\bigcup_{\psi\in\Psi_x}\psi(|\zeta|<R_x).\eqno(2.4)$$
By the proof of Lemma 2.6, there can be no 0-dimensional
components of $\psi_1(|\zeta|<R_x)\cap\psi_2(|\zeta|<R_x)$.  Thus for
$y\in\tilde
V_x$, there is a unique irreducible germ of a variety, $W$ which is
contained in
$\tilde V_x$ and which contains $y$.  Thus there is a Riemann surface $\cR$
and an
injective holomorphic mapping $\chi:\cR\to\tilde V_x$; this Riemann surface
is the
normalization of the singularities of $\tilde V_x$ (see [Ch, \S6]).

\proclaim Proposition 2.7.  Suppose that $f$ does not preserve volume.  If
$x\in
J^*$ has period $n$, and if one of the multipliers of $Df^n_x$ has modulus 1,
then (\dag) does not hold.

\give Proof.  We may assume that $x\in J^*$ is a fixed point of $f$.  If
$\mu_1,\mu_2\in\C$ denote the multipliers of $Df_x$, then we may suppose that
$|\mu_1|<|\mu_2|=1$.  Let us suppose that (\dag) holds at $x$, and let
$\chi:\cR\to\tilde V_x$ be as above.  Set $\tilde x=\chi^{-1}x$.  Then $f$
induces
a biholomorphic mapping $F:=\chi^{-1}\circ f\circ\chi:\cR\to\cR$, and $F(\tilde
x)=x$.  For nonconstant $\psi\in\Psi_x$ we may write
$\psi(\zeta)=x+\sum_{m=k}^\infty a_m\zeta^m$ with $a_m\in\C^2$ and $a_k\ne0$.
Given $\psi'(\zeta)=x+\sum_{m=k'}^\infty a'_m\zeta^m\in\Psi_x$, there exists
$\psi\in\Psi_x$ such that
$$\psi'=\tilde f(\psi)=f\circ\psi(\lambda_\psi^{-1}\zeta)=
\lambda_\psi^{-k} (Df_x\cdot
a_k)\,\zeta^k+O(\zeta^{k+1}).$$
We conclude that $k'=k$, and $a'_k=\lambda_\psi^{-k}(Df_x\cdot a_k)$.  Thus
$a_k$ is an eigenvector of $Df_x$, so $a'_k=\lambda^{-k}_\psi\mu_ja_k$ for
one of
the eigenvalues $\mu_j$.  Since we may choose $\psi$ such that $a_k$ has
maximal length, it follows that $|a_k'|\le|a_k|$,
so $|\mu_j|\ge|\lambda^k_\psi|\ge1$, and so the eigenvalue must be $\mu_2$.

We may also write $\chi(t)=\sum_{m=k}^\infty b_m t^k$ for $b_m\in\C^2$,
$b_k\ne0$.  We compute that $F'(\tilde x)=\mu_2$.  Thus $F:\cR\to\cR$ is an
automorphism with a fixed point $\tilde x$ with multiplier $e^{2\pi i\theta}$.
Passing to covering spaces, we may assume that $\cR$ is $\C$ or the unit
disk.  If
$\theta\in\Q$, then we may assume that $\mu_2=1$, and thus $F$ is the
identity. But this is not possible, since this would mean that $f$ is the
identity
on $\tilde V_x$; but the fixed points of $f$ are discrete.  The other
possibility,
$\theta\notin\Q$ is
also not possible.  For in this case it follows from [BS2, Proposition 2]
that $x$ is
contained in the interior of $K^+$, so $x\notin J^*$.
\qed

We consider the following condition on a point $x\in J^*$:
$$\hbox{\rm Every mapping in $\Psi_x$ is non-constant.}\eqno(\ddag)$$
Note that if $R_x>1$, then $\psi$ is holomorphic on
$\{|\zeta|<R_x\}\supset\{|\zeta|\ge1\}$, and
$\max_{|\zeta|\le 1}G^+\circ\psi=1$.  Since $G^+\circ\psi(0)=0$, it follows
that
$(\ddag)$ holds.  The failure of (\ddag) thus implies
that $R_x=1$ and thus by Proposition 2.2 $\hat\lambda(x,n)=1$ for $n\le0$.
\proclaim Lemma 2.8.  If (\ddag) holds, there
exist $\epsilon>0$ and $0<r<1$ such that for each $\psi\in\Psi_x$ there
exists $\rho\le r$ such that
$$\hbox{\rm dist}(\psi(\zeta),\psi(0))\ge \epsilon$$ for all
$|\zeta|=\rho$.  If (\ddag) holds for all $x\in J^*$, the $\epsilon$ and
$r$ may be
chosen to hold for all $x\in J^*$.

\give Proof.  We expand each $\psi\in\Psi$ in a power series about
$\zeta=0$, $\psi(\zeta)=x+\alpha_1\zeta+\alpha_2\zeta^2+\cdots$,
with $\alpha_j\in\cx2$.  For each $j$, $\psi\mapsto\alpha_j$ is a
continuous mapping from $\Psi$ to $\cx2$.  Since $\psi$ is not constant,
there exist $r_\psi,\epsilon_\psi>0$ such
that $|\psi(\zeta)|>\epsilon_\psi$ for $|\zeta|=r_\psi$.  This inequality
continues to hold in a small neighborhood of $\psi$ inside $\Psi$.  Thus
we obtain $r$ and $\epsilon$ by the compactness of $\Psi$.  \qed

The following shows that if $f$ is quasi-expanding, then each germ $V_x$ is
contained
in a variety
$V({x,\epsilon})$ with uniformly large inner diameter and uniformly bounded
area.
This is an easy consequence of Lemma 2.8 and the fact that ${\rm
Area}(\psi(D))=\int_D|\psi'|^2$.

\proclaim Proposition 2.9.  If $f$ is quasi-expanding, then there exist
$\epsilon>0$
and $A<\infty$ such that for each $x\in J^*$, $V(x,\epsilon)$ is a (closed)
subvariety
of
$B(x,\epsilon)$, the area of $V(x,\epsilon)$ is bounded by $A$.

The following is a strong converse to Proposition 1.5.
\proclaim Theorem 2.10.  If $f$ is not quasi-expanding, then there exists a
point $x\in
J^*$ such that either (\ddag) fails, or $r_x>0$.  In either case, (1.12) fails.

\give Proof.  If $f$ is not quasi-expanding, the by Corollary 2.5 there exists
$\psi\in\Psi$ such that $|\lambda(\psi,n)|=1$ for $n\le0$.  If we set
$\psi_n=\tilde
f^n\psi$ and  $m_{\psi_n}(r)=\sup_{|\zeta|<r}G^+\circ\psi_n(\zeta)$, then
$m_{\psi_n}(|\lambda(\psi,n)|r)=d^nm_\psi(r)$.  Let $\psi_0$ be a limit of
$\psi_{n_j}$ for some subsequence $n_j\to-\infty$.  It follows that
$m_{\psi_0}(1)=0$.  Set $x=\psi_0(0)$.  If $\psi_0$ is constant, then
$(\ddag)$ fails
at  $x$.

Otherwise, if $\psi_0$ is non-constant, we set
$V_0:=\psi_0(\{|\zeta|<1\})$.  Thus
$V_0$ is a connected neighborhood of $x$ in $\tilde V_x\cap\{G^+=0\}$,
where $\tilde
V_x$ is as in (2.4).  Note that $\tilde V_x\cap\{G^+<1\}\subset
J^-\cap\{G^+<1\}$ is
bounded, and thus $\tilde V_x\cap\{G^+<1\}$ is Kobayashi hyperbolic.  Let
$D_K(c)$
denote the disk with center $x$ and radius $c$ in the Kobayashi metric of
$\tilde
V_x\cap\{G^+<1\}$.  We may choose $c$ small enough that $D_K(c)\subset
V_0$.  Now
choose $r$ small enough that the length of $[0,r]$ with respect to the
Kobayashi
metric of the unit disk is less than $c$.  It follows that for any
$\psi\in\Psi_x$
we have $\psi(\{|\zeta|<1\})\subset \tilde V_x\cap\{G^+<1\}$, and thus
$\psi(\{|\zeta|<1\})\subset D_K(c)\subset V_0$.  Thus $r_x\ge r$.
\qed

\section 3.  Area Bounds and Distortion

In this Section we establish a bounded area distortion theorem and use it
to give
sufficient conditions for quasi-expansion.
Recall that if $A\subset\C$ is a doubly connected domain, then $A$ is
conformally equivalent to a circular annulus
$\{\zeta\in\C:r_1<|\zeta|<r_2\}$.  The {\it modulus} of this annulus,
written ${\rm Mod}(A)$, is equal to $\log(r_2/r_1)$.  We will use the
notation $B_{R}=B(0,R)$ for the ball centered at the origin in $\C^n$.

\proclaim Theorem 3.1.  Let $D\subset\C$ be a disk, let $0<R_0<R_1$ be
given, and let
$\phi:D\to B_{R_1}$ be a proper holomorphic map.  Let $A$ denote the
area of the image $\phi(D)$.  The set $\phi^{-1}(B_{R_0})$ is a union of
topological disks.  Let $C$ be any component of $\phi^{-1}(B_{R_0})$.
The set $D-C$ is a topological annulus, and
$${\rm Mod}(D-C)\ge {\log(R_1/R_0)\over {A\over R_1^2}(2 +
{1\over\log(R_1/R_0)} )}.$$

\give Proof.  The modulus of the
annulus $D-C$ is equal to the extremal length of the family of curves
that connect the boundary components.  We recall the computation of
extremal length (see Fuchs [F]).  Given a conformal metric
$\rho(z)|dz|$ on the annulus, the length of a curve $\gamma$ is
$Length(\gamma)=\int_{\gamma}|dz|$.  We define the extremal length
$L$ of the curve family by the formula
$${1\over L}=\inf_{\rho}{{\rm Area}(\rho_0)\over m^2},{\rm\ \
where\ \ }{\rm Area}(\rho_0)=\int\rho^2(z)dA,$$
and $m$ is the infimum of $Length(\gamma)$ for all $\gamma$ in the
curve family.

Any particular choice of $\rho=\rho_0$ gives a lower estimate: $L\ge
m^2/{\rm Area}(\rho_0)$.  Let $ds$ be the Euclidean metric on
$\C^n$.  Let
$r(z)=|z|$ be the radial distance of a point $z\in\C^n$ to the origin,
and let $\rho_0$ be the pullback under $\phi$ of the metric which is
defined as
$ds/r$ on $B_{R_1}-B_{R_0}$ and which is zero on
$B_{R_0}$.

We will estimate the minimal length of a curve
and the area for the metric $\rho_0$.
Define $g(v)=\log|v|={1\over2}\log(v,v)$ for $v\in\C^n$, so that
$\nabla g={v/|v|^2}$.  Let $\gamma(t)$ be a path in $D-C$ with
$\gamma(0)\in\partial C$ and
$\gamma(1)\in\partial D$.
$$\eqalign {
{\rm Length}(\gamma) = & \int_0^1{|(\phi\circ\gamma)'(t)|
\over|(\phi\circ\gamma)(t)|} dt
=\int_0^1{|(\phi\circ\gamma)'(t)|\cdot |(\phi\circ\gamma)(t)|
\over|(\phi\circ\gamma)(t)|^2} dt \cr
&\ge\int_0^1(\phi\circ\gamma)'\cdot\nabla g\, dt=g(1)-g(0)\cr
&=\log|\phi\gamma(1)| -\log|\phi\gamma(0)|=\log(R_1/R_0).\cr}$$

Let $F(r)$ denote the area of $\phi(D)\cap B_r$ with respect to the
standard metric on $\C^n$.  By definition, $F(R_1)=A$.  The area of
$D-C$ with respect to $\rho_0$ is
$${\rm Area}(\rho_0) =\int_{R_0}^{R_1} {F'(r)\over r^2}dr.$$
Now we integrate by parts and use the property (see [Ch, p.\ 189]) that
$F(r)/r^2$ is non-decreasing in $r$ to obtain
$$\eqalign{
\int_{R_0}^{R_1}{F'(r)\over r^2} dr = \left({F(R_1)\over
R_1^2} -{F(R_0)\over R_0^2}\right) -\int_{R_0}^{R_1} F(r){-2\over
r^2}dr\cr
={F(R_1)\over R_1^2} - {F(R_0)\over R_0^2} + 2\int_{R_0}^{R_1}{F(r)\over
r^2}{dr\over r}\cr
\le {F(R_1)\over R_1^2} + {2F(R_1)\over R_1^2}\log{R_1\over R_0}.}$$
Thus
$$\eqalign{
{\rm Mod}(D-C)=L\ge{ (\log{R_1\over R_0})^2\over {A\over R_1^2} +
{2A\over R_1^2}\log{R_1\over R_0}}\cr
={\log{R_1\over R_0}\over {A\over R_1^2}(2+{1\over\log{R_1\over
R_0}})}\cr }$$
which is the desired estimate.  \qed

Theorem 3.1 yields the following {\it Bounded Area Distortion Theorem,} which
is of independent interest.  In dimension $n=1$, the fact that the area is
weighted by
the multiplicity of the mapping becomes crucial.  For
$n=1$, the part of Theorem 3.2 concerning the containment
$\{|\zeta|<ar\}\subset D_0\subset\{|\zeta|<r\}$ is given in [CJY, Lemma 2.2].
\proclaim Theorem 3.2.  Let $A<\infty$ and $\chi>1$ be given.  Then
there exist $\rho>0$ and $a>0$ with the following property:  If
$\cD\subset\C$ is a simply connected domain containing the origin, and if
$\phi:\cD\to B_R$ is a proper holomorphic mapping with $\phi(0)=0$ and
${\rm Area}(\phi(\cD))\le A$, then for some $r$ the component $D_0$ of
$\phi^{-1}(B_{\rho R})$ containing the origin satisfies
$$\{|\zeta|<ar\}\subset D_0\subset\{|\zeta|< r\}\subset\{|\zeta|<\chi
r\}\subset\cD.$$

\give Proof.  We define $r>0$ to be the minimum value for
which
$D\subset\{|\zeta|<r\}$, and we define $\kappa>0$ to be the maximum value for
which $\{|\zeta|<\kappa r\}\subset\cD$.  There are points $\zeta_1\in\partial
D$, $|\zeta_1|=r$ and $\zeta_2\in\partial\cD$, $|\zeta_2|=\kappa r$.  By
Teichm\"uller's Theorem [A, Theorem 4-7], the modulus of the annulus $\cD-\bar
D$ is no greater than the modulus of the complement in $\C$ of the segments
$[-r,0]$ and $[\kappa r,+\infty)$.  By Theorem 3.1 we have
$${\rm Mod}(\C-([-1,0]\cup[\kappa,+\infty))) \ge {\log (R/\rho)\over {A\over
R^2}(2+{1\over \log(R/\rho)})}.$$
The quantity on the left hand side of this inequality has been much studied.
The estimate given by equation (4-21) of [A] shows that
the modulus is less than ${(2\pi)}^{-1}\log(16(\kappa+1))$.  Thus if we take
$\rho$ sufficiently small then we can obtain
$\kappa\ge\chi$.

Now define $a>0$ as the largest value for which $\{|\zeta|<ar\}\subset D$.
Then there is a point
$\zeta_0\in\partial D$ with $|\zeta_0|=ar$.  Without loss of generality, we
may assume that $\zeta_0>0$.  By the Cauchy estimate, we have
$|\phi'(\zeta)|\le R/((\kappa-1)r)$ for $|\zeta|\le r$.  Since $\phi(0)=0$ and
$\phi(\zeta_0)\in\partial B_{\rho R}$, we have
$$\rho R\le \int_0^{\zeta_0}|\phi'(t)|\, dt\le {R\over
(\kappa-1)r}|\zeta_0|={Ra\over \kappa-1}.$$
We conclude that $a\ge (\kappa-1)\rho$, which completes the proof.
\qed

We consider the condition:
$$\eqalign{ {\rm There\ exist\ }\epsilon>0{\rm\ and\ } A<\infty {\rm\
such\ that\ }
{\rm for\ all\ }x\in S,\cr
  V({x,\epsilon}){\rm\ is\ closed\ in\
}B(x,\epsilon),{\rm\ \ and\ \ \  }{\rm
Area}(V({x,\epsilon}))<A.\cr}\eqno(3.1)$$
If (3.1) holds, then by Bishop's Theorem (see [Ch, p.\ 205]) the family
$\{V({x,\epsilon}):x\in S\}$ is pre-compact in the Hausdorff topology.  By
Lemma 2.8,
there is a unique extension to a family of varieties
$\{V({x,\epsilon}):x\in  J^*\}$, with $V({x,\epsilon})\subset J^-$.
We also consider the related condition, in which
$S$ is replaced by
$ J^*$:
$$\eqalign{ {\rm There\ exist\ }\epsilon>0{\rm\ and\ } A<\infty {\rm\
such\ that\ }
{\rm for\ all\ }x\in  J^*,\cr
  V({x,\epsilon}){\rm\ is\ closed\ in\
}B(x,\epsilon),{\rm\ \ and\ \ \  }{\rm
Area}(V({x,\epsilon}))<A.\cr}\eqno(3.1')$$
Clearly $(3.1')\Rightarrow (3.1)$.  By Proposition 2.9, conditions (3.1)
and $(3.1')$
hold if $f$ is quasi-expanding.

\proclaim Theorem 3.3.  $(3.1)\Rightarrow$ ((\ddag) holds for all
$x\in J^*$) $\Rightarrow(3.1')$.

\give Proof. Suppose that (3.1) holds.  For $\chi>1$, let $\rho$ and $a$ be the
constants from Theorem 3.2.  We may assume that
$\epsilon>0$ is small enough that $\max_{\bar B(x,\epsilon)}G^+\le 1$ for
all $x\in J^*$.  For $x\in S$ we let $\psi_x\in\psi_S$ be the corresponding
immersion.  Let $\cD$ (respectively,  $D_0$) denote the connected component
of $\psi_x^{-1}(B(x,\epsilon))$ (respectively,
$\psi_x^{-1}(B(x,\rho\epsilon))$) containing the origin.  Since
$G^+\circ\psi_x\le 1$ on $\{|\zeta|<r\chi\}\subset\cD$ it follows that
$r\chi\le1$, so $r\le1/\chi$.  We note that
$dist(\psi_x(\zeta),x)=dist(\psi_x(\zeta),\psi_x(0))=\rho\epsilon$ for all
$\zeta\in\partial D_0$, and that $\partial D_0\subset\{|\zeta|\le1/\chi\}$.
It follows that if a sequence of mappings $\psi_{x_j}$ converges to a map
$\psi$, then there will be a point $\zeta$ with $|\zeta|\le1/\chi$ and
$dist(\psi(\zeta),\psi(0))=\rho\epsilon$, and so $\psi$ cannot be constant.
Thus
(\ddag) holds.

Next let us suppose that (\ddag) holds.  We have
already observed that with $\epsilon$ as in Lemma 3.3, $V({x,\epsilon})$ is
closed in $B(x,\epsilon)$.  By the compactness of
$\Psi$, $C_0=\sup_\psi\sup_{|\zeta|\le r}|\psi'(\zeta)|<\infty$, so
$$\hbox{\rm Area}(V({x,\epsilon}))\le\int_{|\zeta|\le
r}|\psi'(\zeta)|^2\le \pi r^2C_0^2,$$
so $(3.1')$ holds.  \qed

The germ $V_x$ being contained in $K$ is equivalent to
$G^+\circ\psi(\zeta)=0$ for
$|\zeta|<r$; by Proposition 1.5 this prevents quasi-expansion.  The
condition that no
germ $V_x$ is contained in
$K^+=\{G^+=0\}$ may be viewed as a weak form of transversality between
$V_x$ and $K^+$.

\proclaim Theorem 3.4.  If  (3.1) holds, and if no germ $V_x$, $x\in J^*$,
is contained in $K$, then $f$ is quasi-expanding.

\give Proof.  Let $\epsilon>0$ and $A<\infty$ be as in (3.1). For $\chi>1$,
let $\rho$ be the constant in Theorem 3.2, corresponding to the number $A$.  By
the continuity of
$G^+$, we may shrink $\epsilon$ so that
$$\max_{x\in J^*}\max_{\bar B(x,\epsilon)}G^+\le 1.$$
We set
$$c_0:=\min_{x\in J^*}\max_{V({x,\rho\epsilon})}G^+.$$
We claim that $c_0>0$.  By Theorem 3.3, $\{V_x, x\in J^*\}$ is a compact
family of
varieties.  If $c_0=0$, then by compactness we would have $G^+=0$ on
$V({x,\rho\epsilon})$ for some
$x\in J^*$.  This germ $V_x$ would be contained in $K$, contradicting our
hypothesis.  Thus $c_0>0$, and we choose $n$ such that $d^nc_0\ge1$.

For $x\in S$, let $\cD$ (respectively $D_0$) be the connected component of
$\psi^{-1}B(x,\epsilon)$ (respectively $\psi^{-1}B(x,\rho\epsilon)$)
containing the origin.  Thus
$$D_0\subset\{|\zeta|<r\}\subset\{|\zeta|<r\chi\}\subset\cD.$$
Since $G^+\circ\psi_x\le 1$ on $\cD$, it follows that $r\chi\le1$.  We also
have
$$c_0\le\max_{\zeta\in  D_0}G^+\le\max_{|\zeta|\le r}G^+.$$
Now applying $f^n$, the set $V_x$ is mapped to $V_{f^nx}$, which is
unformized by $\psi_{x_n}$.  It follows that
$$\max_{|\zeta|\le|\lambda(x,n)|r}G^+\circ\psi_{x_n}\ge c_0d^n\ge1.$$
By the normalization (1.3) on $\psi_{x_n}$ it follows that
$|\lambda(x,n)|r\ge1$.  We conclude that
$$|\lambda(x,n)|r\ge1\ge\chi r,$$
or $|\lambda(x,n)|\ge\chi>1$.  By Theorem 1.2, then, $f^n$ is
quasi-expanding.  Thus by Proposition 1.3, $f$ is quasi-expanding.
\qed

The following two Corollaries are just restatements of Proposition 2.3 and
Theorem
2.8 in terms of the families
$\psi_S$ given as examples at the end of \S1.

\proclaim Corollary 3.5.  A necessary and sufficient condition for $f$
to be quasi-expanding is that there exist $\epsilon>0$ and $A<\infty$
such that for each $\delta>0$ there is an $\eta>0$
such that for each saddle point $p$ we have:
$W^u_\epsilon(p)$ is closed in
$B(p,\epsilon)$, ${\rm Area}(W^u_\epsilon(p))\le A$, and
$\sup_{W^u_\delta(p)}G^+\ge\eta$.

\proclaim Corollary 3.6.  A necessary  and sufficient condition for $f$
to be quasi-expanding is that there exist $\epsilon>0$, $A<\infty$
and saddle points $p,q$ such that for each $\delta>0$ there is an
$\eta>0$ such that
$W^u_\epsilon(z)$ is closed in
$B(z,\epsilon)$,
${\rm Area}(W^u_\epsilon(z))\le A$, and $\sup_{W^u_\delta(z)}G^+\ge\eta$
for all $z\in
W^u(p)\cap W^s(q)$.

\proclaim Corollary 3.7.  The condition that $f$ is quasi-expanding is
independent of the family $\psi_S$.

A quasi-expanding mapping has a
certain uniform contraction along backward orbits, at finite scale.

\proclaim Theorem 3.8.  If $f$ is quasi-expanding, then there exist
$\epsilon>0$ and $\theta<1$ such that for all $x\in J^*$ and $n\ge0$,
${\rm diam}(f^{-n}V({x,\epsilon}))\le\theta^n$ and ${\rm
Area}(f^{-n}V({x,\epsilon}))\le\theta^n$.

\give Proof.  Let $A<\infty$ and $\epsilon>0$ be as in Proposition
2.9.  Let $c$ be the maximum of $G^+$ in an
$\epsilon$-neighborhood of $ J^*$.  Let $M(r)$ be as in Theorem 1.2, and
let $r_c$ be such that $M(r_c)=c$.  We note that by choosing $\epsilon$
sufficiently small, we can make $c$ arbitrarily close to zero.  By the
continuity of $M(r)$, then, the constant $r_c$ may be taken arbitrarily
close to zero.

Let  $a$ and $\rho$ be as in Theorem 3.2.  For $x\in J^*$, let
$\cD_x$ denote the connected component of $\psi_x^{-1}V({x,\epsilon})$
containing the
origin.  It follows from Theorem 3.2 that there exists $r_x$ such that
$$\{|\zeta|<ar_x\}\subset \psi_x^{-1}V({x,\rho\epsilon})
\subset\{|\zeta|< r_x\}.$$
It is evident that $ar_x\le r_c$.  Define
$$M:=\sup_{\psi\in\Psi}\sup_{|\zeta|<a^{-1} r_c}|\psi'(\zeta)|<\infty$$
which is finite by the compactness of $\Psi$.
By the transformation formula (1.4), and by (4) of Theorem 1.2, we have
$$f^{-n}V({x,\rho\epsilon})\subset f^{-n}\psi_x(|\zeta|<a^{-1} r_c)
\subset\psi_{f^{-n}x}(|\zeta|<\kappa^{-n}a^{-1} r_c).$$
Thus we obtain the estimates
$$\eqalign{ {\rm Area}(f^{-n}V({x,\rho\epsilon}))\le
\int_{|\zeta|<\kappa^{-n}a^{-1} r_c}|\psi'_{f^{-n}x}(\zeta)|^2\cr
\le \pi(\kappa^{-n} a^{-1} r_c)^2 M^2.\cr}$$
and
$$\eqalign{ {\rm diam}_{\C^2}(f^{-n}V({x,\rho\epsilon}))\le (\max|\psi'|)
{\rm diam}\{|\zeta|<\kappa^{-n}a^{-1} r_c\}\cr
\le M\cdot 2\kappa^{-n}a^{-1} r_c.\cr}$$
Finally, it suffices to take $r_c$ sufficiently small
that $2Ma^{-1} r_c<1$. \qed
\proclaim Proposition 3.9.  If $f$ is uniformly hyperbolic on $J^*$, then
$f$ and $f^{-1}$ are quasi-expanding.

\give Proof.  Let $\cW^u=\{W^u(x):x\in J^*\}$ denote the lamination defined
by the unstable manifolds through points of $J^*$.  Since $\cW^u$ contains
the sets $V(x,\epsilon)$, it follows that $(3.1')$ holds.  For $x\in J^*$
let $\Vert\cdot\Vert_x^\#$ denote the metric on the tangent space $E^u_x$, as
well as the distance induced on $W^u(x)$ by $\Vert\cdot\Vert_x^\#$.  For $x\in
J^*$, let $0\le r(x)\le\infty$ denote the largest number such that the
$\Vert\cdot\Vert^\#_x$-disk $D(x,r(x))\subset W^u(x)$ is contained in
$\{G^+=0\}$.  Since $|\lambda_\psi|\ge1$, we have $r(fx)\ge r(x)$.

By [BS7, \S5], $x\mapsto\Vert\cdot\Vert^\#_x$ varies continuously.  It follows
that $J^*\ni x\mapsto r(x)$ is upper semicontinuous.  We will show that
$r(x)=0$.  Suppose, to the contrary, that $R:=\sup_{x\in J^*}r(x)>0$.  By
upper semicontinuity there exists $x_0\in J^*$ with $r(x_0)=R$.  If
$R=\infty$, it follows that $W^u(x_0)$ is conformally equivalent to $\C$, and
that $W^u(x_0)\subset K$, which is a contradiction.  If $R<\infty$, then we
let $X$ denote the $\omega$-limit set of $x_0$, i.e., the set of limits of
sequences $\{f^{n_j}x_0\}$, $n_j\to\infty$.  Clearly $X$ is
$f$-invariant and compact, so it is a hyperbolic set for $f$.  Since $r$ is
upper semicontinuous, we have $r(x)=R$ for all $x\in X$.  Thus for
$x_0\in X$
$$W^u(x_0)\supset \bigcup_{n\ge
0}f^n(D(f^{-n}x_0,r(f^{-n}x_0)))=\bigcup_{n\ge0}f^n(D(f^{-n}x_0,R)).$$
Since $\Vert\cdot\Vert^\#$ is comparable to euclidean distance, there exists
$r_0>0$ such that $D(f^{-n}x_0,R)$ contains a disk with Euclidean radius $r_0$
in $W^u(f^{-n}x_0)$.  It follows, then, as in [BS1] that
$W^u(x_0)\subset\{G^+=0\}$ and is conformally equivalent to $\C$, which is a
contradiction.  Thus we conclude that $r=0$ on $J^*$.  It now follows from
Theorem 3.5 that $f$ is quasi-expanding.  The argument for $f^{-1}$ is the
same.
\qed

\proclaim Proposition 3.10.  Suppose $f$ is quasi-expanding.  Then there exist
$\epsilon>0$ and $N<\infty$ such that $f^nV({x,\epsilon})\supset
V({f^nx,\epsilon})$ for $x\in J^*$ and $n\ge N$.

\give Proof.   By Lemma 3.3 and Theorem 3.2 there exist $\epsilon>0$ and
$L<\infty$ such that
$$\psi_x(\{\zeta\in\C:|\zeta|<1/L\})\subset V({x,\epsilon})\subset
\psi_x(\{\zeta\in\C:|\zeta|<L\})$$
for each $\psi_x\in\Psi_x$.  Let $\kappa>1$ be as in Theorem 1.2, and choose
$N$ such that $\kappa^N>L^2$.  It follows, then that $N$ has the desired
property.
\qed

\section 4.  Equivalence of Families of Metrics

In this Section we will show that two families of mappings are
quasi-expanding.  These are the real mappings of maximal entropy
(Theorem 4.6) and the (uniformly) hyperbolic mappings (Corollary 4.11),
which gives an alternate proof of Proposition 3.9.  We will establish
quasi-expansion in both cases by showing that the metric
$\Vert\cdot\Vert^\#$ is uniformly expanded by $f$ (condition 4 of Theorem
1.2).  The metric
$\Vert\cdot\Vert^\#$ is in general not equivalent to the Euclidean metric.
We present another metric
$\Vert\cdot\Vert^{(L)}$ which, like $\Vert\cdot\Vert^\#$, is defined only in
terms of the complex structure and the function $G^+$.  The metric
$\Vert\cdot\Vert^{(L)}$ is uniformly expanded by real mappings of maximal
entropy.  We develop tools to compare the expansion with respect to
different metrics, specifically $\Vert\cdot\Vert^{(L)}$ and the Euclidean
metric, and use this to show that certain mappings are also expanding with
respect to the metric
$\Vert\cdot\Vert^\#$.

We continue to let $E^u_x\subset T_x\C^2$ denote the unstable tangent
space.  We will compare several norms on $E^u_x$. By $\|\cdot\|^e$ we
denote the norm on $E^u_x$ induced by the euclidean metric on $\C^2$.  If
$\Vert\cdot\Vert$ denotes an complex affine invariant metric on
$E^u_x$, then $\Vert\cdot\Vert$ is determined by its value at the origin.
Since
$E^u_x$ has complex dimension 1, it follows that all such affine metrics
are real multiples of each other, i.e.\ there exists $\alpha=\alpha(x)$
such that $\Vert\cdot\Vert=\alpha\Vert\cdot\Vert^e$.

Now we define another metric.  For $0<L<\infty$, we let $D^{(L)}=D^{(L)}_x$
denote the connected component
of $\{\zeta\in\cx{}:G^+\circ\psi_x(\zeta)<L\}$ which contains the origin.
Since $G^+\circ\psi_x$ is subharmonic on $\cx{}$, it follows from the
maximum principle that $D^{(L)}$ is simply connected.  We let
$ds_P$ denote the Poincar\'e metric of $D^L$ at $\zeta=0$, and we define
$\Vert v\Vert_x^{(L)}=ds_P(v )$ for any tangent vector $v$. If $L\le L'$,
then $D^{(L)}\subset D^{(L')}$, and so the Poincar\'e metrics of the two
domains have the reverse inequality, so that $\Vert\cdot\Vert^{(L)}\ge
\Vert\cdot\Vert^{(L')}$. By the identity $G^+\circ f=d\cdot G^+$ it follows
that
$f$ maps the set
$D^{(L)}_x$ to the set $D^{(dL)}_{fx}$. Thus for $v\in E^u_x$ we have
$$\Vert v\Vert_x^{(L)}=\Vert Df_xv\Vert_{fx}^{(dL)}. \eqno(4.1)$$

If $\Vert\cdot\Vert^1_x$ and
$\vert\cdot\Vert^2_x$ are families of affine metrics for $x\in S$, there is a
comparison function $a^1_2: S\to\R$ defined by
$a^1_2(x)=\log(\Vert v\Vert^1_x/\Vert v\Vert_x^2)$, where $v$ is any nonzero
element of $E^u_x$.  We say that $\Vert\cdot\Vert^1$ is {\it equivalent}
to
$\Vert\cdot\Vert^2$ if $a^1_2(x)$ is a bounded function of $x\in S$.  Note that
this is an equivalence relation.

Let us observe that $\Vert\cdot\Vert^\#$ is
equivalent to $\Vert\cdot\Vert^{(1)}$.  To see this, we recall from
the definition of $\Vert\cdot\Vert^\#_x$ that the unit disk in
$\cx{}$ is the largest disk centered at the origin which is contained in
$\{G^+\circ\psi_x<1\}$, i.e. $D^{(1)}$ contains the disk
$\{\zeta\in\cx{}:|\zeta|<1\}$, and there is a point
$\zeta_0\in\partial D^{(1)}$ with $|\zeta_0|=1$.  Since $G^+\circ\psi_x$ is
subharmonic on $\cx{}$, it follows that $D^{(1)}$ is simply connected.  Now let
$\chi:\{|\zeta|<1\}\to D^{(1)}$ be a conformal equivalence such that
$\chi(0)=0$.  It follows that the Poincar\'e metric $ds_P$ satisfies:
$ds_P({\bf
1})=|\chi'(0)|$.  By the Koebe $1/4$-Theorem, we have ${1\over
4}\le|\chi'(0)|\le
1$.  It follows from the definition of $\Vert\cdot\Vert^\#$ that
$${1\over 4}\Vert\cdot\Vert_x^\#\le
\Vert\cdot\Vert_x^{(1)}\le\Vert\cdot\Vert_x^\#.\eqno(4.2)$$

We may compare the metrics $\Vert\cdot\Vert^{(L)}$ and
$\Vert\cdot\Vert^{(L')}$,
$L<L'$ as follows.  We define conformal maps $\chi_x^{(L)}:D_x^{(L)}\to B_1$,
where $B_1=\{|z|<1\}$, and $\chi_x^{(L)}(0)=0$.  Since $L<L'$, we may
define the
induced map
$$\rho:=\chi_x^{(L')}\circ(\chi_x^{(L)})^{-1}:B_1\to B_1\eqno(4.3)$$
satisfies
$$a_{L}^{L'}:=\log{\Vert v\Vert^{(L)}_x \over \Vert v\Vert
^{(L')}_x}= -\log|\rho'(0)|.\eqno(4.4)$$
{}From the fact that $D_x^{(L)}$ is a strict subset of $D_x^{(L')}$, it
follows that $|\rho'(0)|<1$, and thus $a_{(L)}^{(L')}<0$, which corresponds
to the earlier observation that
$\Vert\cdot\Vert^{(L)}_x<\Vert\cdot\Vert^{(L')}_x$ for all $x$.  A lower
bound on $a_{(L)}^{(L')}$ is equivalent to a lower bound on $\rho'(0)$,
which is equivalent to the existence of $0<r<1$ such that
$\{|z|<r\}\subset\rho(B_1)$ for all $x\in S$.

\proclaim Lemma 4.1.  There exists $0<L_0<\infty$, (depending on $f$), such
that
for $L_0<L'<\infty$,
$\Vert\cdot\Vert^{(L_0)}$ is equivalent to
$\Vert\cdot\Vert^{(L')}$.

\give Proof. We let
$\pi$ denote the coordinate projection onto the second coordinate axis.  It
follows that the restrictions of
$|\pi|$ and
$G^+$ to $J^-$ are proper exhaustions of $J^-$.  Let us choose $C_1$ such
that
$$\{G^+=0\}\cap J^-\subset\{|\pi|<C_1\}.$$
For $L_0$ sufficiently large, there exists $C_2>C_1$ such that
$$\{G^+>L\}\subset\{|\pi|>C_2\}.$$
Finally, for $L'<\infty$ we may choose $C_3>C_2$ sufficiently large that
$$\{G^+<L'\}\subset \{|\pi|<C_3\}.$$

Let $s$ denote the distance between $\{|z|=C_1\}$ and $\{|z|=C_2\}$, measured
with respect to the Poincar\'e metric on $\{|z|<C_3\}$.  Let $0<r<1$ be
chosen so
that the Poincar\'e distance from $0$ to $r$ inside $B_1$ is equal to $s$.  We
will show that $\rho(B_1)$ contains the disk of radius $r$ about the origin.
For if $|z_0|<r$, then the Poincar\'e distance between $0$ to $z_0$ is less
than
$s$.  We consider $A:=\pi\circ\psi_x\circ(\chi^{(L')})^{-1}(0)$ and
$B:=\pi\circ\psi_x\circ(\chi^{(L')})^{-1}(z_0)$.  Since
$0\in J\subset J^-\cap\{G^+=0\}$, it follows that $|A|<C_1$.  By the definition
of $D^{(L')}$ and $C_3$, it follows that
$$\pi\circ\psi_x\circ(\chi^{(L')})^{-1}(B_1)\subset\{|z|<C_3\}.$$
Thus $\pi\circ\psi_x\circ(\chi^{(L')})^{-1}$, as a mapping from $B_1$ to the
disk $\{|z|<C_3\}$ decreases the respective Poincar\'e metrics.  Thus the
distance between the points $A$ and $B$ is less than $s$, so we conclude
that $B$
is contained in the disk $\{|z|<C_2\}$.  By the definition of $C_2$, then, it
follows that $G^+(\psi_x\circ(\chi^{(L')})^{-1}(z_0))<L$.  Thus we conclude
that
$z_0$ is in the range of $\rho$, which gives the desired lower estimate.
\qed

Given a family of metrics $\Vert\cdot\Vert$ the effect of the differential $Df$
is measured by the function $c(x,n)=\log(\Vert Df^n_xv\Vert_{fx}/\Vert
v\Vert_x)$, where $v$ is a nonzero element of $E^u_x$.  The chain rule
gives the
following cocycle condition for $c$:
$$c(x,n+m)=c(x,n)+c(f^nx,m).$$
\give Remark.  We observe that the ``natural'' metrics $\Vert\cdot\Vert^\#$ and
$\Vert\cdot\Vert^{(L)}$ are expanded under $f$, although the expansion
is not uniform in $x\in S$.  If $c^\#=c^\#(\cdot,1)$ denotes the cocycle
corresponding to the metric $\Vert\cdot\Vert^\#$, then we
have $c^\#(x,1)=\log|\lambda_x|$, so $c^\#>0$.

For the metric $\Vert\cdot\Vert^{(L)}$ we note that by (4.2) the corresponding
cocycle satisfies
$$c^{(L)}(x,1)=\log{\Vert Df_xv\Vert_{fx}^{(L)}\over \Vert v\Vert_x^{(L)} }
= \log{\Vert v \Vert^{(L/d)}_x\over \Vert v\Vert_x^{(L)} } =
a_{L}^{L/d}(x).\eqno(4.5)$$
By (4.4) we have $c^{(L)}>0$.

Given two families of metrics $\Vert\cdot\Vert^1$ and $\Vert\cdot\Vert^2$, the
corresponding cocyles are related by the coboundary equation:
$$c^1(x,1)-c^2(x,1)=a^1_2(fx)-a^1_2(x).$$
We say that two cocyles $c^1$ and $c^2$ are {\it equivalent} if they
satisfy the coboundary equation
$$c^1(x,1)-c^2(x,1)=\alpha(fx)-\alpha(x)$$
for some function $\alpha: S\to\R$ which is bounded.  With these definitions,
equivalent families of metrics produce equivalent cocycles. Indeed the above
equation can be solved with the bounded function $a^1_2$ for $\alpha$.

\proclaim Lemma 4.2.  For all $0<L,L'<\infty$, the cocycle $c^L$ is
bounded, and
$\Vert\cdot\Vert^{(L)}$ is equivalent to $\Vert\cdot\Vert^{(L')}$.

\give Proof.  For a point $x_0\in S$ and $v_0\in E^u_{x_0}$, we set
$x_j=f^jx_0$ and $v_j=Df^jx_0$.  Applying (4.1), we have
$$\Vert Dfv_n\Vert_{x_{n+1}}^{(L)} = \Vert v_n\Vert_{x_n}^{(d^{-1}L)} = \Vert
Dfv_{n-1}\Vert_{x_n}^{(d^{-1}L)}=\cdots =\Vert Dfv_0\Vert_{x_1}^{(d^{-n}L)}.$$
This gives
$$c^{(L)}(x_{n+1},1) = \log{ \Vert Dfv_n\Vert_{x_{n+1}}^{(L)}\over\Vert
v_n\Vert_{x_n}^{(L)}} = \cdots = \log{ \Vert Df v_0\Vert_{x_1}^{(d^{-n}L)}
\over
\Vert v_0\Vert_{x_0}^{(d^{-n}L)} } = c^{(d^{-n}L)}(x_0,1).$$

By Lemma 4.1 and the monotonicity of $\Vert\cdot\Vert^{(L)}$, there exists
$\kappa<\infty$ such that $0\le a_{L_1}^{L_0}\le \kappa$ for all $L_0\le L_1\le
L_2$. This gives
$0\le c^{(L_1)}(x,1)\le\kappa$ for all $x\in  S$.  Thus $0\le
c^{(d^{-n}L_1)}(x,1)\le\kappa$ for all $x\in  S$.

Now choose $n$ such that $L<L'\le d^nL$.  It follows that
$$0\le a_{L'}^L\le a_{d^nL}^L = a_{d^nL}^{d^{n-1}L} +
a_{d^{n-1}L}^{d^{n-2}L} +\cdots +a_{dL}^L = c^{(d^nL)} +\cdots+
c^{(dL)}\le n\kappa,$$
which gives the equivalence between the metrics.
\qed

\proclaim Lemma 4.3.  A cocycle which is equivalent to a bounded cocycle
is a bounded cocycle.

\give Proof.  If $c^1$ and $c^2$ are equivalent, then
$$c^1(x,1)-c^2(x,1)=\alpha(fx)-\alpha(x)$$
for some bounded function $\alpha$.  If $c^2$ is bounded, then so is $c^1$,
since
each term on the right hand side of
$$c^1(x,1)=a^1_2(fx)-a^1_2(x)+c^2(x,1)$$
is bounded.\qed

\proclaim Corollary 4.4.  The cocycle $c^\#$ is bounded, i.e. there exists
$\chi<\infty$ such that $|\lambda_x|<\chi$ for all $x\in S$.

\give Proof.  By Lemma 4.2, the cocycle $c^{(1)}$ is bounded.  And by (4.2)
$\Vert\cdot\Vert^\#$ is equivalent to $\Vert\cdot\Vert^{(1)}$.  By Lemma 4.3 it
follows that
$c^\#(x,1)=\log|\lambda_x|$ is bounded. \qed

A cocycle $c$ is said to be {\it eventually positive} if for some $n>0$ and
some $K>0$ we have $c(x,n)\ge K$ for all $x\in S$.
We first observe:

\proclaim Lemma 4.5.  A cocycle which is boundedly cohomologous to an
eventually positive cocycle is eventually positive.

\give Proof.  If $\Vert\cdot\Vert^1$ is eventually positive, then $c^1(x,n)\ge
K$ for $n\ge k$.  If $\Vert\cdot\Vert^2$ is equivalent, then
$$c^1(x,1)-c^2(x,1)=a(fx)-a(x)$$
for some bounded comparison function $a$.  Now
$$c^1(x,n)-c^2(x,n)=a(f^nx)-a(x)$$
where $a$ is bounded, say $|a|\le C$.  If $nK<2C$, then
$$c^1(x,nk)=c^2(x,nk)=a(f^{nk}x)-a(x)$$
so
$$c^2(x,nk)=-c^1(x,nk)+a(f^{nk}x)-a(x)\ge nK -2C>0.$$
\qed

A cocycle is {\it immediately positive} if $c(x,1)\ge K>0$.  We recall
that one of the equivalent conditions in the definition of
quasi-expanding is that the cocycle $c^\#$ corresponding to the metric
$\Vert\cdot\Vert^\#$ is immediately positive.  By Proposition 1.5, it
follows that $c^\#$ is eventually positive if and only if it is immediately
positive.

Now we consider mappings $f$ which are real.  This means that the real
subspace ${\bf R}^2$ is invariant under $f$, or in terms of coordinates,
$f$ commutes with complex conjugation, i.e., $f(\bar x,\bar y)=\bar
f(x,y)$.  We
let $f_R$ denote the restriction of $f$ to ${\bf R}^2$.  We will say that
$f$ is
a real mapping with maximal entropy if the real restriction $f_{{\bf R}^2}$ has
entropy equal to $\log d$.  Several results from [BLS] apply to real mappings
with maximal entropy. In this case it follows that $J\subset {\bf R}^2$, that
$J= J^*$, and the periodic points are dense in $J$.  Thus, if $p$ is a
(real) periodic point, we may further normalize the uniformizing mapping
$\psi_p:\cx{}\to W^u(p)$ so that $\psi_p({\bf R})\subset{\bf R}^2$.  In this
case, it follows that $\psi_p(\cx{})\cap J\subset{\bf R}^2$, and thus
$\psi_p^{-1}(J)\subset {\bf R}$.

\proclaim Theorem 4.6.  If $f$ is a real mapping of maximal entropy, then the
cocyle corresponding to the metric
$\Vert\cdot\Vert^{(L)}$ satisfies $c^{(L)}(x,1)\ge \log d$.  Further, $f$
and
$f^{-1}$ are both quasi-expanding.

\give Proof. We observed above that if $f$ is a real mapping of maximal
entropy, then for each saddle point $\psi_x^{-1}(J)\subset\R$.  Thus by
Proposition 4.7,  $c^{(L)}$ is a positive cocyle. Since
$\Vert\cdot\Vert^{(L)}$ is equivalent to $\Vert\cdot\Vert^\#$ it follows from
Lemma 4.5 that $c^\#$ is eventually positive.  By (4) of Theorem 1.2, some
iterate $f^N$ is quasi-expanding.  By Proposition 1.3, then,
$f$ itself is quasi-expanding.  The argument for $f^{-1}$ is similar.
\qed

\proclaim Proposition 4.7.  If $x\in S$ is such that
$\psi^{-1}_xJ$ is contained in a straight line in $\C$, then
$c^{(L)}(x,1)\ge\log d>1$.  If, in addition,
$J\cap W^u(x)$ is not connected, then we have $c^{(L)}(x,1)>\log d$.

\give Proof.  Without loss of generality, we may assume that the line is
$\R\subset\C$.  We will estimate $c^{(L)}$ as in (4.5).  To do this, we let
$h^+$ denote the unique continuous function on $D^{(L)}$ with the following
properties:  $0\le h^+<L$, $h^+=0$ on $D^{(L)}\cap{\bf R}$,
$h^+$ is harmonic on $D^{(L)}-{\bf R}$, and $h^+$ takes the boundary limit $L$
at all points of $(\partial D^{(L)})-{\bf R}$.  Since
$\psi^{-1}_xJ\subset{\bf R}$, it follows from the maximum principle that
$h^+\le G^+\circ\psi_x$ on $D^{(L)}$.  Thus $D^+:=\{h^+<L/d\}\supset
D^{(L/d)}$.  Let $\chi^+:D^+\to B_1$ denote the conformal mapping such that
$\chi^+(0)=0$ and $\chi^+(0)'>0$.  If we set $\rho^+:=\chi^+\circ
(\chi^{(L)})^{-1}:B_1\to B_1$, then as in (4.4), we have the estimate
$$c^{(L)}=-\log|\rho'(0)|\ge- \log|\rho^+(0)'|.$$

We will show that $\rho^+(0)'=1/d$, which gives $c^{(L)}\ge \log d$.  In
order to
do this, we let $H^+$ be the function on $B_1$ which is the image of $h^+$
under $\rho^+$.  We note that since $f$ is real, the set $D^{(L)}$ is
invariant under complex conjugation, and thus $\rho^+$ commutes with
conjugation.  Thus the real function $H^+$ is invariant under complex
conjugation.  Further, $H^+$ has the properties of being equal to zero on
the axis $(-1,1)$, harmonic on $B_1-(-1,1)$, and taking boundary values
$L$ on the non-real points of $\partial B_1$.  Let $\varphi$ denote the
conformal mapping from $B_1$ to the strip $\{\zeta\in\C:-1<\Im( \zeta)<1\}$
such that the upper/lower portion of $\partial B_1$ is taken to
$\{\Im(\zeta)=\pm 1\}$.  It follows that $H^+=L|\Im(\varphi)|$.  The image of
$\rho^+$ is given by
$\rho^+(B_1)=\{H^+<L/d\}$.  Thus $\rho^+$ is given by $\varphi^{-1}\circ
g_d\circ\varphi$, where $g_d(z)=z/d$.  We may assume that
$\varphi(0)=0$, so it follows from the fact that $g_d'=1/d$ that
$\rho^+(0)'=1/d$.

If $J\cap W^u(x)$ is not connected, then $h^+<G^+\circ\psi_x$ on
$D^{(L)}$ because
$\{h^+=0\}=\psi_x^{-1}(J)\ne\{G^+\circ\psi_x=0\}=\R$.    Thus
$c^{(L)}=-\log|\rho'(0)|>- \log|\rho^+(0)'|=d$.
\qed
Next we give an improved estimate for one-sided points, which play an important
role in [BS].
\proclaim Proposition 4.8.  If $x\in S$ and $\psi^{-1}_xJ$
is contained in a half-line, then $c^{(L)}(x,1)\ge 2\log d$.  If, in addition,
$\psi^{-1}_xJ$ is not connected, then $c^{(L)}(x,1)>2\log d$.

\give Proof.  Without loss of generality, we may assume that $\psi^{-1}_xJ$
contains the origin and is contained in the positive half-line $[0,\infty)$.
The proof now proceeds along the lines of the proof of Proposition 4.7 with
the modification that the function $|\Im(\zeta)|$ is replaced by
$|\Im(\sqrt\zeta)|$.  Let $S(t)=\{\zeta\in\C:|\Im(\sqrt\zeta)|<t\}$.  We let
$\varphi$ denote the conformal mapping from
$B_1$ to the set $S(1)$.  The factor of 2 enters
because $d^2$ is the multiplier which maps
$S(L/d)$ to $S(L)$.
\qed

An affine metric on $\cx{}_x$ induces a distance function on $W^u(x)$ via the
mapping $\psi_x:\cx{}_x\to W^u(x)$.  The metric
$\Vert\cdot\Vert^\#$ induces the distance
$dist^\#(\psi_x(\zeta_1),\psi_x(\zeta_2))=|\zeta_1-\zeta_2|$.  Any other metric
is of the form $\Vert\cdot\Vert'=a\,\Vert\cdot\Vert^\#$, and the induced
distance
is given as $dist'(\psi_x(\zeta_1),\psi_x(\zeta_2))=a(x)|\zeta_1-\zeta_2|$.
Given a metric, we let $\Delta_x\subset W^u(x)$ denote the unit disk in
$W^u(x)$ with center at $x$.  We say that the metric is {\it admissible}
if there are constants $0<c'<c''<\infty$ such that the diameter, measured
with respect to the Euclidean  metric on $\cx2$ satisfies
$$c'\le diam_{\cx2}(\Delta_x)\le c''$$
for all $x\in S$.

Admissibility of a metric is not a strictly local property since it
involves the
immersions $\psi_x$.  If we work with the metric
$\Vert\cdot\Vert^\#$, then the boundary
$\Delta_x$ contains a point of $\{G^+=1\}\cap J^-$.  Since this is a
compact set, we have an upper bound on the diameter of $\Delta_x$.  The
lower bound on the diameter follows because this set is at positive
distance from $J$ (which contains
$x$).  Thus
$\Vert\cdot\Vert^\#$ is admissible.

It need not be true that a metric equivalent to an admissible metric is
itself admissible.  For $0<\tau\le 1$, let us consider the scaled metric
$\tau\Vert\cdot\Vert^\#$.  With respect to this metric, $\Delta_x$
satisfies $\psi_x^{-1}\Delta_x=\{\zeta\in\cx{}_x:|\zeta|<{1\over\tau}\}$.
It is evident that $\Psi$ is a normal family if and only if the functions
$\psi_x$ are bounded on this set for each $\tau>0$.  For fixed $\tau$, this
gives an upper bound on $\sup_{x\in S}diam(\Delta_x)$.  (The case $\tau=1$
is already a lower bound.)  Thus we see that $f$ is quasi-expanding if and
only if $\tau\Vert\cdot\Vert^\#$ is admissible for every $0<\tau<1$.
In other words, if $f$ is not quasi-expanding, then $\tau\|\cdot\|^\#$ is
not admissible for some $0<\tau<1$.

\proclaim Lemma 4.9.  Any two admissible metrics are equivalent.

\give Proof. Suppose that $\Vert\cdot\Vert$ and $\Vert\cdot\Vert'$ are
admissible metrics.  If they are not equivalent, we may choose a sequence $x_k$
such that $\Vert\cdot\Vert_{x_k}=\epsilon_k\Vert\cdot\Vert'_{x_k}$, and
$\lim_{k\to\infty}\epsilon_k=0$.  Let $D_k:=\{\zeta\in\cx{}_{x_k}:\Vert
\zeta\Vert_{x_k}<1\}$.  Thus $D_k':=\{\zeta\in\cx{}_{x_k}:\Vert
\zeta\Vert'_{x_k}<1\}=\epsilon_k D_k\subset D_k$.  Let $\delta_k$ denote the
diameter of
$D_k'$, measured with respect to the Kobayashi metric of $D_k$.

Since $\Vert\cdot\Vert$ is admissible, there is a bounded set $B\subset\cx2$
such that $\Delta_k=\psi_{x_k}(D_k)\subset B$ for all $k$.  Since the Kobayashi
metric decreases under holomorphic mappings,  the diameter of
$\psi_{x_k}(D_k')$,
measured with respect to the Kobayashi metric of $B$ is no larger than
$\delta_k$.  Further, since $B$ is bounded, the Kobayashi metric of $B$
dominates the Euclidean metric.  Thus for some constant $C<\infty$, the
Euclidean diameter of $\Delta_k'=\psi_{x_k}(D_k')$ is no larger than $\delta_k
C$.  But if $\epsilon_k\to0$, it follows that the relative diameter $\delta_k$
also tends to zero. Thus the Euclidean diameters of $\Delta_k'$ are not bounded
below, which contradicts the admissibility of $\Vert\cdot\Vert'$.  This
contradiction shows that the two metrics must be equivalent. \qed

\proclaim Theorem 4.10.  If $f$ is uniformly hyperbolic on $J^*$, then
$\Vert\cdot\Vert^e$ is an admissible metric.

\give Proof.  Let $\cW^u$ denote the lamination of $W^u(J^*)$ by unstable
manifolds.  Each unstable manifold is uniformized by ${\bf C}$, and thus has a
unique complex affine structure.  It was shown in [BS7] that this affine
structure varies continuously.  For each $p\in J$, we may assign  an affine
metric on $W^u(p)$ by using the metric $\Vert\cdot\Vert_p^e$, induced by
$\phi_p:{\bf C}\to W^u(p)$.  By the continuity of the affine structure, the
sets $\Delta_p=\psi_p\{\zeta\in\cx{}:\Vert\zeta\Vert_p^e<1\}$ vary
continuously.  In particular, their diameters will be bounded above and
below in terms of the euclidean metric on $\cx2$.
\qed
We conclude with another proof of Proposition 3.9.
\proclaim Corollary 4.11.  If $f$ is uniformly hyperbolic on $J^*$, then $f$
and $f^{-1}$ are quasi-expanding.

\give Proof. Let $c^e$ denote the cocycle corresponding to the Euclidean
metric.  If $f$ is uniformly hyperbolic, then $c^e$ is eventually positive.
Further, since both $\Vert\cdot\Vert^\#$ and  (by Theorem 4.10)
$\Vert\cdot\Vert^e$ are admissible, they are equivalent by Lemma 4.9.
By Lemma 4.5, the cocycle $c^\#$ is eventually positive.  By (4) of
Theorem 1.2 and  Proposition 1.3 it follows that $f$ is quasi-expanding.  \qed

\section 5.  Local Folding

In this Section we show how conditions (\dag) and (\ddag) express
themselves in terms
of local folding.  In \S2 we showed how (\ddag) corresponds to a bound on
the local
area of the varieties $\cV$.  Here we show (Propositions 5.1--3) how it
corresponds to
a bound on the local folding of $\cV$.

For $\psi\in\Psi_x$, we define ${\rm
Ord}(\psi)=\min\{n\ge1:\psi^{(n)}(0)\ne0\}$.  Thus
${\rm Ord}(\psi)<\infty$ if and only if $\psi$ is non-constant.  If $j={\rm
Ord}(\psi)<\infty$, then $\psi(\zeta)=x+a_j\zeta^j+\cdots$, where we set
$a_j=\psi^{(j)}(0)/j!$.  By $E_x$ we denote the complex linear span of $a_j$ in
$T_x\C^2$.  $E_x$ coincides with the tangent cone of the variety $V_x$ at
$x$ (see [Ch
\S8]). By Lemma 2.6, $E_x$ is independent of the choice of $\psi$.
In the following discussion of folding, we will use the notation $E_x$ to
denote the
complex affine line passing through $x$ in the direction of the tangent
cone of $V_x$ at
$x$.  Let $\pi:\C^2\to E_x$ denote a complex affine projection map.  Let
$\psi\in\Psi_x$ be non-constant.  For an open set $\cN\subset E_x$, we let
$\omega$
denote the connected component of $\psi^{-1}(\pi^{-1}\cN)$ containing $x$.
Since $E_x$
is the tangent cone to $V(\psi)$ at $x$, we may choose $\cN$ sufficiently
small that
$\omega$ is relatively compact inside $\Delta$. For $\hat\psi\in\Psi$, we let
$V(\hat\psi,\cN)$ denote the connected component of
$\psi(\Delta)\cap\pi^{-1}\cN$
containing $\psi(0)$.  If $\hat\psi\in\Psi$ is uniformly close to $\psi$ in a
neighborhood of $\bar\omega$, then $V(\hat\psi,\cN)$ is a subvariety of
$\pi^{-1}\cN$, and $\pi|V(\hat\psi,\cN):V(\hat\psi,\cN)\to\cN$ is proper.  If
$y=\hat\psi(0)$, then in analogy with \S2, we may write $V(y,\cN)$ for
$V(\hat\psi,\cN)$.

Let us define $\tau(x):=\sup_{\psi\in\Psi_x}{\rm Ord}(\psi)$.  If
$\tau(x)=\infty$ then
(\ddag) does not hold.

\proclaim Proposition 5.1.  Suppose (\dag) holds at $x$ and
$\tau(x)=\infty$.  Then
for each $k<\infty$ and for an arbitrarily small neighborhood $\cN$ of $x$
inside $E_x$,
there exists $y\in J^*\cap\pi^{-1}\cN$ such that $V(y,\cN)$ is a
nonsingular subvariety
of $\pi^{-1}\cN$, and
$\pi|V(y,\cN):V(y,\cN)\to\cN$ is proper with mapping degree $\ge k$.

\give Proof.  Since $\tau(x)=\infty$, $\Psi_x$ contains elements with
arbitrarily
high order.  Thus for each $k$, there exists $\psi\in\Psi_x$ with $m:={\rm
Ord}(\psi)\ge k$.  Let us choose $\{p_j\}\subset S$ such that
$\lim_{j\to\infty}\phi_{p_j}=\psi$.  As was observed above, we may choose
$j$ large
and $\cN$ sufficiently small that $\pi|V({p_j},\cN):V({p_j},\cN)\to\cN$ is
proper.  Since $\phi_{p_j}\in\psi_S$, the varieties $V({p_j},\cN)$ are
regular.
The map
$\zeta\mapsto\pi\circ\psi(\zeta)$ is $m$-to-1 near $\zeta=0$.  It follows
that $m$ is
the mapping degree of $\pi|V({p_j}):V({p_j})\to\cN$. \qed
\proclaim Proposition 5.2.  Suppose that $x$ satisfies (\dag), and suppose
that for
each sufficiently small neighborhood $\cN$ of $x$ inside $E_x$ there exists
$y$ close to
$x$ such that $\pi|V(y,\cN):V(y,\cN)\to\cN$ is a proper map of degree $k$.
Then
$\tau(x)\ge k$.  Further $\tau(x)$ is the smallest number with this property.

\give Proof.  We choose a sequence of neighborhoods $\cN_j$ decreasing to
$\{x\}$ and
let $\psi_j$ be the corresponding functions.  Passing to a subsequence, we
may suppose
that $\psi_j\to\psi\in\Psi_x$.  For $\epsilon>0$, we may choose $\cN$ small
enough that
$\psi^{-1}(\pi^{-1}\bar\cN)\subset\{|\zeta|<\epsilon\}$.  Thus for $j$
large enough we
have $\phi_j^{-1}(\pi^{-1}\bar\cN)\subset\{|\zeta|<\epsilon\}$.  It follows
that
$\pi\circ\phi_j:\{|\zeta|<\epsilon\}\to\cN$ is a $k$-to-one mapping.  Thus
$\psi:\{|\zeta|<\epsilon\}\to\cN_\epsilon$ is $k$-to-one.  Since this holds
for all
$\epsilon>0$, it follows that ${\rm Ord}(\psi)=k$.

To establish the final statement, we suppose first that $\tau(x)=\infty$.
Then by
Proposition 5.1, there are $\psi\in\Psi_x$ yielding branched covers of
degree $\ge k$.
If $\tau(x)<\infty$, we may choose $\psi\in\Psi_x$ with ${\rm
Ord}(\psi)=\tau(x)$. Again, by the argument of Proposition 5.1, there is a
local
branching of order
$k=\tau(x)$.\qed

For a positive integer $k$, we set $\cJ_k=\{x\in J^*:\tau(x)=k\}$.  Thus
$\cJ_1,\cJ_2,\dots$ is a partition of $\{x\in J^*:\tau(x)<\infty\}$.  Since
$J^*\ni
x\mapsto \tau(x)$ is upper semicontinuous, the set $\bigcup_{k\ge m}\cJ_k$ is
closed (and $\bigcup_{k<m}\cJ_k$ is open) in $\{x\in J^*:\tau(x)<\infty\}$
for each
$m$.

Figure 1 illustrates the case where $\tau(x)=k$, and $\psi_x(\C)$ is a
nonsingular manifold.  By Lemma 2.6, there is a neighborhood $\cN$ of $x$
inside $V_x$ ($\pi^{-1}\cN$ is shaded in  Figure 1) with the following
properties.  If $y\in\cJ$ is sufficiently close to $x$ and
$\psi\in\Psi_y$, then the variety $V_\psi(\cN)$ is a $j$-fold branched
cover over $\cN$.  Since we are working in the complex domain a manifold
$\psi(\C)$ cannot lie ``to one side" of $\psi_x(\C)$.  Also highlighted is
a regular point where $\psi(\C)$ has a vertical tangent.  On compact sets
outside the shaded neighborhood, $\psi(\C)$ has the geometry of $j$
distinct manifolds which approach $\psi_x(\C)$ in the $C^1$ topology as
$y\to x$.

We may interpret the mapping degree of $\pi|V_\psi(\cN)$ as measuring
the local folding of the variety $V_\psi(\cN)$ at $x$.  The following
result asserts that the maximal amount of local folding at $x$ is given by
$\tau(x)$,
which also measures the maximal order of  vanishing of the derivatives of
the parametrizations.

\epsfxsize3.4in
\centerline{\epsfbox{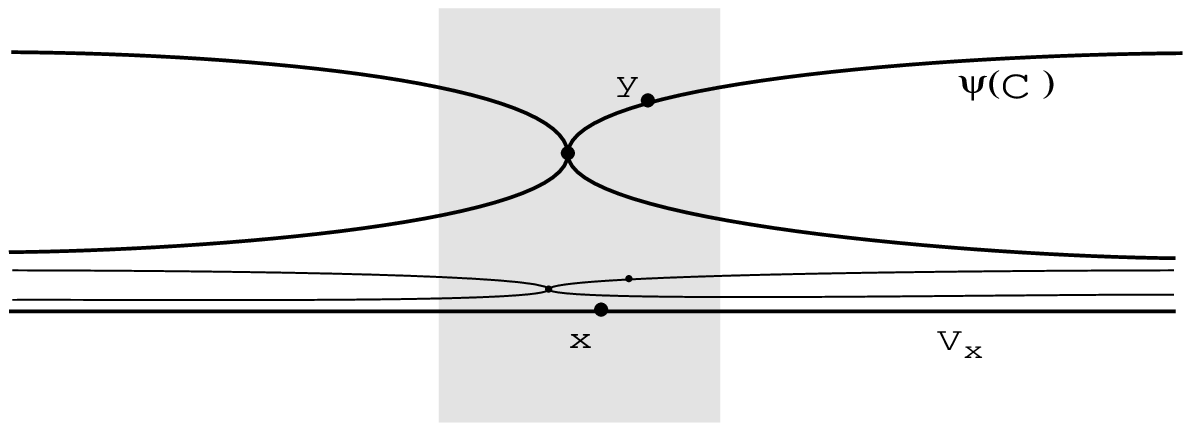}}
\centerline{Figure 1}

\give Proposition 5.3.  If $x\in\cJ_1$, there are neighborhoods
$x\in\cU_0\subset\cU\subset\C^2$ such that
$\{V(y,\epsilon)\cap\cU:y\in\bar\cU_0\cap
J^*\}$ is a lamination.  If $x\in J^*-\cJ_1$, then there is no such
lamination at $x$.

\give Proof.   For
$x\in\cJ_1$, every $\psi\in\Psi_x$ has nonvanishing differential at
$\zeta=0$.  Since $\Psi$ is a normal family, it follows that the images
$V(x,\epsilon)$ are regular and form a lamination.  Conversely, if $x_0\in
J^*-\cJ_1$, then for any open sets $\cU_0\subset\cU$ containing $x_0$, there
will be varieties $V(x,\epsilon)$, $x\in \cU_0\cap J^*$ which project as in
Figure 1.  Thus there is no lamination at $x_0$.  \qed

\section 6.  Expansion

For $x\in J^*$ and $k<\infty$, let us define $\Psi^k_x=\{\psi\in\Psi_x:{\rm
Ord}(\psi)=k\}$, and
$$\gamma_k(x):=\sup_{\psi\in\Psi^k_x}\left|{\psi^{(k)}(0)\over k!}\right|,$$
which is finite by the normality of $\Psi$.  We have
$\cJ_k\subset\{\gamma_k>0\}$, and the set of all
points of $J^*$ where (\dag) holds coincides with
$\bigcup_{k\ge1}\{\gamma_k>0\}$.  By
Lemma 2.6 and the normality of
$\Psi$, it follows that
$\{\gamma_k>0\}\ni x\mapsto E_x$ is continuous.  For $x\in J^*$ with
$\gamma_k(x)>0$,
we define
$$\Vert v\Vert^{\#,k}_x:={|v|/\gamma_k(x)}{\rm\ \ \ for\ }v\in E_x.$$
Since $\Psi$ is generated by the normal limits of elements of $\psi_S$, we
have
$$\Vert\
\Vert^{\#,1}_x=\liminf_{p\in S, p\to x}\Vert\ \Vert^{\#}_p,$$
where $\Vert\ \Vert^\#_p$ was defined in \S1.  Since $\gamma_k(x)$ is upper
semicontinuous, we have an upper bound $m_k:=\sup_{x\in
J^*}\gamma_k(x)<\infty$, so we
have a lower bound in terms of the euclidean metric:
$${|v|\over m_k}\le\Vert v\Vert^{\#,k}_x,{\rm\ \ \ for\ }v\in E_x.$$
If
$\phi\in \Psi^k_x$ (resp.\
$\phi_1\in\Psi^k_{fx}$) realizes the supremum defining $\gamma_k$ at $x$
(resp.\ $fx$),
then
$$|\lambda_{\tilde f^{-1}\phi_1}|^k\le\Vert Df^n_x\Vert^{\#,k}:={ \Vert
Df^n_xv\Vert^{\#,k}_{f^nx} \over \Vert
v\Vert^{\#,k}_x}\le|\lambda_\phi|^k.\eqno(6.1)$$

If (\ddag) holds, then $x\in\cJ_k$ for $k=\tau(x)$, and we define a metric
$\Vert\
\Vert^\#_x$ on $E_x$ by setting $\Vert v\Vert^\#_x:=\Vert v\Vert^{\#,k}_x$
for $v\in
E_x$.  This replaces the definition given in \S1; in general the two
definitions may
disagree for  $x=p\in S$.  If $S$ is a compact subset of $\cJ_k$, then
$c:=\inf_{x\in S}\gamma_k(x)>0$.  Thus for any compact $S\subset\cJ_k$
$$m_k^{-1}|v|\le\Vert v\Vert^\#_x\le c^{-1}|v|, {\rm\ for\ }x\in S, v\in
E_x\eqno(6.2)$$ gives an equivalence between $\Vert\ \Vert^\#_x$ and the
euclidean
metric.
\proclaim Proposition 6.1.  If $\gamma_k(x)>0$, then $\Vert
Df^n_x|E_x\Vert^e\le
C\gamma(x,n)$ for $n\le0$.  If $n_j\to\infty$ is a sequence with
$f^{n_j}x\to\hat
x\in\cJ_k$, then there exists $c>0$ such that $\Vert Df^{n_j}_x|E_x\Vert^e\ge
c\lambda(x,n_j)$.

\give Proof.  By the definition of $\Vert\ \Vert^{\#,k}$, we have
$$\Vert Df^n_x|E_x\Vert^e=\lambda(x,n){\gamma_k(f^nx)\over\gamma_k(x)}.$$
Thus $C=m/\gamma_k(x)$ is our desired bound.  If $\hat x\in\cJ_k$, then
$\eta:=\liminf_{x\in\cJ_k,x\to\hat x}\gamma_k(x)>0$.  Thus if ${n_j}$ is
sufficiently
large, we have $\gamma_k(f^{n_j}x)>\eta/2$, which gives the desired estimate
with $c=\eta(2\gamma_k(x))^{-1}$. \qed

\give Remark on Expansion.  We may interpret the Proposition as follows.
$Df^n_x|E_x$ has uniform contraction along the backward orbit of a point
$x\in\cJ_k$.
If there is a sequence of times $n_j\to\infty$ such that
$dist(f^{n_j}x,\{\tau(x)<k\})$
is bounded below, then $Df^n|E_{f^nx}$ has exponential growth during the
times $n=n_j$.

\proclaim Theorem 6.2.  Let $f$ be quasi-expanding, and let $\nu$ be an
ergodic invariant measure supported on
$\cJ$.  Then the Lyapunov exponent of $\nu$ satisfies
$\Lambda(\nu)\ge\log\kappa>0$.
If $\cJ_k$ has full measure for $\nu$, then $\Lambda(\nu)\ge k\log\kappa>0$.

\give Proof.  The Lyapunov exponent of the measure $\nu$ is given by the
formula
$$\Lambda(\nu)=\lim_{n\to\infty}{1\over n} \int\log\Vert Df^n_x\Vert\,\nu(x).$$
Since the family $E^u_x$ is invariant, it follows that
$${1\over n}\log\Vert Df^n_x\Vert \ge{1\over n}\Vert Df^n_x|E^u_x\Vert
={1\over n}\sum_{j=0}^{n-1}\log\Vert Df_{f^jx}|E^u_{f^jx}\Vert.$$
By the invariance of $\nu$, we have $\int\log\Vert
Df_{f^jx}|E^u_{f^jx}\Vert\nu(x)=\int\log\Vert Df_{x}|E^u_{x}\Vert\nu(x)$, so
$${1\over n}\int\log\Vert Df^n_x\Vert\,\nu(x) \ge{1\over n}\int\log\Vert
Df^n_x|E^u_x\Vert\,\nu(x)=\int\log\Vert
Df_x|E^u_x\Vert\,\nu(x).$$

It will suffice to consider the case when all the mass of $\nu$ is on
$\cJ_k$.  For $x\in\cJ_k$ we have $\gamma_k(x)>0$, so
$$\Vert Df_x|E^u_x\Vert ={|Df_x v|_{fx}\over |v|_x} = \Vert Df_x\Vert^\#\
{\gamma_k(fx)\over\gamma_k(x)}.$$
By [LS, Proposition 2.2], we have
$$\int\log{\gamma_k(fx)\over\gamma_k(x)}\,\nu(x)=0$$
It follows from (6.1) that
$$\int\log\Vert Df_x|E^u_x\Vert\,\nu(x)=\int\log\Vert
Df_x|E^u_x\Vert^\#\,\nu(x)\ge\log\kappa^k,$$
and the last inequality follows from (6.1).
\qed

\proclaim Corollary 6.3. If $f$ is quasi-expanding and quasi-contracting,
then every
periodic point in $J^*$ is a saddle point.  Further, there is a $\kappa>1$
such that if
$\lambda^+$ and
$\lambda^-$ denote the larger and smaller eigenvalues of $Df^n$ at a saddle
point of
period $n$, then $|\lambda^-|\le\kappa^{-n}<\kappa^n\le|\lambda^+|$.

Let us use the notation
$\cJ_k'=\{x\in\cJ_k:\alpha(x)\cap\cJ_k\ne\emptyset\}$, where
$\alpha(x)$ is the $\alpha$-limit set, i.e.\ the accumulation points of
sequences
$f^{n_j}x$ with $n_j\to-\infty$.  By the Poincar\'e Recurrence Theorem,
$\cJ_k'$ has full measure for any invariant measure on $\cJ_k$.  If
$\cJ_k$ is compact (which occurs, for instance, if $k=\sup_{x\in
J^*}\tau(x)<\infty$), then $\cJ'_k=\cJ_k$.

\epsfxsize3in
\centerline{\epsfbox{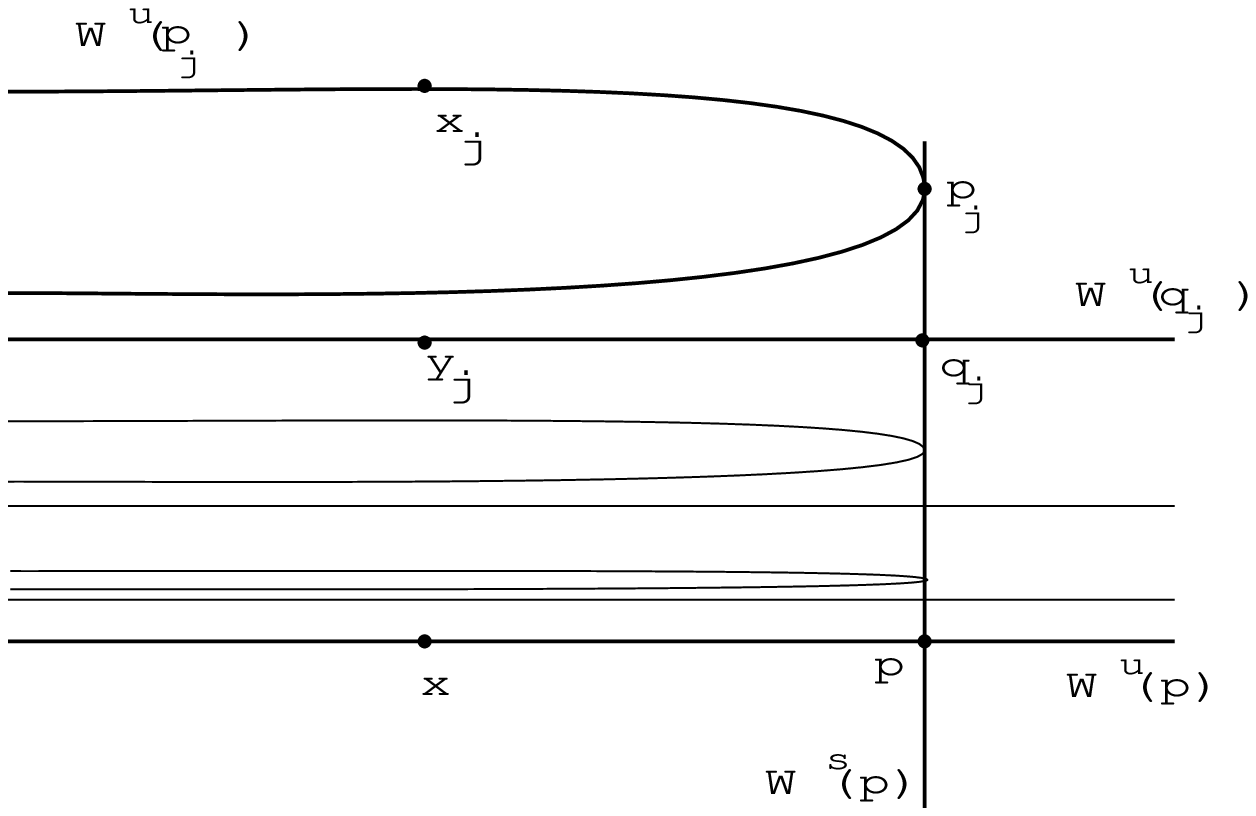}}
\centerline{Figure 2.}

\proclaim Lemma 6.4.  Let $x\in\cJ_k'$ be a point of backward contraction.  If
$\psi\in\Psi_x$, then the number of critical points of $\psi$, counted with
multiplicity, is no greater than $k-1$; and for $\zeta\in\C$, the number of
preimages
$\psi^{-1}(\psi(\zeta))$, of a point $\psi(\zeta)$ is no greater than $k$.

\give Proof.  Let $\zeta_1,\dots,\zeta_j$ be critical points of $\psi$.
Then $\tilde f^{-n}(\psi)$ has critical points at
$\lambda(x,-n)\zeta_i$ for $1\le i\le j$.  Let
$\tilde\psi\in\Psi_y$ denote any normal limit of a subsequence of
$\tilde f^{-n}(\psi)$.  By the backward contraction, all the critical
points converge
to the origin in the limit, so it follows that
${\rm Ord}(\tilde
\psi)$ is one greater than the sum of the multiplicities of the critical
points
$\zeta_1,\dots,\zeta_j$.  Since $\alpha(x)\cap\cJ_k\ne\emptyset$, we
may take the subsequence such that $y\in\cJ_k$.  Thus ${\rm
Ord}(\tilde\psi)\le k$, and thus the total orders of the critical points
must be less than $k$.  A similar argument shows that
$\#\psi^{-1}(\psi(\zeta))\le k$.  \qed

In Figure 2 we suppose that $p$ is a saddle point and that the unstable
manifolds $W^u(p_j)$ and $W^u(q_j)$ are smooth.  Let $\phi_p:\C\to
W^u(p)\subset\C^2$, $\phi_{p_j}:\C\to W^u(p_j)\subset\C^2$, and
$\phi_{q_j}:\C\to W^u(q_j)\subset\C^2$ be holomorphic embeddings.  Let
$\beta_j,\gamma_j\in\C$ be such that $\phi_{p_j}(\beta_j)=x_j$,
$\phi_{q_j}(\gamma_j)=y_j$, and $\phi_p(\beta)=x$.  It follows that
$\phi_{q_j}\to\phi_p\in\Psi_p$, $\gamma_j\to\beta$, and
$\phi_{q_j}(\zeta+\gamma_j)\to\phi_p(\zeta+\beta)\in\Psi_x$.  If
$W^u(p_j)$ and $W^u(p)$ have a simple (quadratic) tangency at $p_j$, then
$\phi_{p_j}(\zeta)\to\phi_p(\alpha\zeta^2)\in\Psi_p$ for some $|\alpha|=1$,
and $\phi_{x_j}(\zeta)\to\phi_p(\alpha(a\zeta+b)^2)\in\Psi_x$, where
$a,b\in\C$ are chosen so that (1.3) holds.  Figure 2 is
consistent with the properties $x\in\cJ_1$ and $p\in\cJ_2$.

We say that $\psi\in\Psi$ is a {\it homogeneous parametrization} if it has
the form $\psi(\zeta)=\phi(c\zeta^k)$, where $c\in\C$, and $\phi:\C\to\C^2$
is an
immersion.  If
$\psi_1,\psi_2\in\Psi^k_x$ are two homogeneous parametrizations,  then by
(1.3) and Lemma 2.6, they differ only by a rotation of the variable $\zeta$.

\proclaim Lemma 6.5.  Let $x\in\cJ_k'$ be a point of backward contraction. If
$\phi:\C\to\C^2$ is an immersion which induces $V_x$, then there is a
polynomial
$p(\zeta)$ of degree no greater than $k$ such that $p(0)=0$ and
$\psi(\zeta)=\phi(p(\zeta))$.  In particular, every $\psi\in\Psi_x^k$ is a
homogeneous
parametrization.

\give Proof.  If $\phi$ is an immersion, then
we claim that $p:=\phi^{-1}\circ\psi:\C\to\C$ is a well-defined
holomorphic mapping.  It is evident that $p$ is analytic on the
domain $\cD:=\{\zeta\in\C:\psi(\zeta)\in\phi(\C)\}$.  Note that we must
have $p(\zeta)\to\infty$ as $\zeta\to\partial\cD$.  Otherwise, if
$\zeta_j\to\zeta_0\in\partial\cD$ and $p(\zeta_j)\to c_0\in\C$ we have
$p(\zeta_0)=\phi^{-1}(\psi(\zeta_0))=c_0$, so $\psi(\zeta_0)=\phi(c_0)$,
which means
that $\zeta_0\in\cD$.  It follows that $1/p$ may be extended to a
continuous function
on $\C$ by setting it equal to $0$ on $\C-\cD$.  By Rado's Theorem, then,
it follows
that $1/p$ is holomorphic on $\C$, which means that
$\C-\cD$ is a discrete set. By the argument above, $p$ has a pole at each
$\zeta_0\in\C-\cD$.  But if
$\lim_{\zeta\to\zeta_0}p(\zeta)=\infty$, then
$\lim_{\zeta\to\infty}\phi(\zeta)=\psi(\zeta_0)$, forcing $\phi$ to be
constant.  This completes the proof of the claim.   By Lemma 6.4, $\psi$ is
at most $k$-to-one, it follows that $p$ is a polynomial of degree no
greater than $k$.

If, in addition, ${\rm Ord}(\psi)=k$, then the multiplicity of the
critical point of the origin is already $k-1$, so we must have
$p(\zeta)=\alpha\zeta^k$.
\qed

\proclaim Theorem 6.6.  Let $\psi_S$ be as in Example 2 in \S1.  Then the
definition of $\Vert\
\Vert^\#_p$ as given in \S1 coincides with the definition given in \S4.

\give Proof.  Let $k$ be such that $p\in\cJ_k'$.  Let $\phi:\C\to  W^u(p)$
denote
the normalized uniformization.  If $\phi(\zeta)=p+a\zeta+\cdots$, then $\Vert
v\Vert^\#_p=|v/a|$, according to the definition in \S1.

If $\psi\in\Psi_p$, it follows from Lemma 6.5 that
$\psi(\zeta)=\phi(c\zeta^k)$ for some scalar $c$ with $|c|=1$.  Since any
two
homogeneous parametrizations agree up to a rotation of parameter,
$a$ must be maximal, so $\Vert v\Vert^\#_p=|v/a|$ according to definition
in \S4. \qed

\proclaim Theorem 6.7.  If $f$ is quasi-expanding, then $\cJ_1$ is an open,
dense subset of $\cJ$.

\give Proof.  Suppose that $k$ is the minimum value of $\tau$ on $\cJ$.
It follows that $\cJ_k=\{\tau<k+1\}$ is an open set.  Since $\cJ_k$ is
$f$-invariant, and since $\cJ_k$ is a nontrivial open subset of the
support of $\mu$, it follows that $\cJ_k$ has full $\mu$ measure and is
thus dense in $\cJ={\rm supp}(\mu)$.
It will suffice to show that $k=1$.

First we claim that for $x\in\cJ_k$, each $g\in\Psi_x$ with ${\rm
Ord}(g)=k$ has the form $g=\psi(\alpha\zeta^k)$.  Let
$x\in\cJ_k$ be a periodic point, and let
$\phi_x:\C\to W^u(x)$ denote the uniformization of the unstable
manifold, normalized to satisfy (1.3).  Let $g_x\in\Psi_x$ be a map
such that ${\rm Ord}(g_x)=k$.  It follows from Lemma 6.5 that
$g_x(\zeta)=\phi_x(\alpha_x\zeta^k)$ with $|\alpha_x|=1$.  For general
$x_0\in\cJ_k$, we may let $x_j$ be a sequence of periodic points
converging to $x_0$.  Passing to a subsequence if necessary, we have
that $g_{x_j}=\phi_{x_j}(\alpha_{x_j}\zeta^k)$ converges to
$g_{x_0}\in\Psi_{x_0}$, and $\alpha_{x_j}\to\alpha$.  It follows that
$\phi_{x_j}$ converges to a function $\phi:\C\to\C$ with
$g_{x_0}(\zeta)=\phi(\alpha\zeta^k)$.

For $x\in S$, let $\psi_x\in\phi_S$.  We know that
$\omega(x)\subset\cJ$, so $\omega(x)\cap\cJ_i\ne\emptyset$ for some
$1\le i\le m$.  By Lemma 6.5 there is a polynomial
$p_x(\zeta)$ of degree no greater than $m$, with $p(0)=0$, such that
$\phi(\zeta)=\psi_x(p_x(\zeta))$.  Thus $g_x=\psi_x(p_x(\alpha\zeta^k))$.

Now let $h\in\Psi_x$ be an element with ${\rm Ord}(h)=k$.  There exist
immersions $\psi_{x_i}\in\psi_S$ which converge to $h$.  For each $i$,
let $g_{x_i}\in\Psi_{x_i}$ be an element such that ${\rm
Ord}(g_{x_i})=k$.  Then as above we have a mapping $\phi_i$ such that
$g_{x_i}(\zeta)=\phi_i(\alpha_i\zeta^k)$.  Since $\phi_i$ and
$\psi_{x_i}$ both have the normalization (1.3), we have
$g_{x_k}=\psi_{x_i}(p_{x_i}(\alpha_i\zeta^k))$.  Since $\Psi$ is a
normal family, we may extract a subsequence so that $g_{x_i}\to
G\in\Psi_x$.

Next we claim that the polynomials $\{p_{x_i}\}$ form a normal family.
Since the degree of $p_{x_i}$ is bounded by $m$, it suffices to show
that $\max_{|\zeta|\le1}|p_{x_i}|$ is bounded.  For each $i$, let
$r_i$ denote the radius of the largest disk centered at the origin and
contained in the image $p_{x_i}(|\zeta|<1)$. Since $\phi_i$ and
$\psi_{x_i}$ are normalized according to (1.3), we must have
$r_i\le1$.  We suppose that $C_i\to\infty$ and derive a contradiction.
We have $\max_{|\zeta|\le1}|p_{x_i}(\zeta)|=1$, so we may extract a
subsequence such that $C_{i}^{-1}p_{x_i}\to q$, a polynomial of degree
no greater than
$m$.   Again we have $q(0)=0$ and $\max_{|\zeta|\le1}|q(\zeta)|=1$, so
that $q$ is non-constant.  Thus $q(|\zeta|\le1)$ contains
neighborhood of the origin.  On the other hand, the interior radius
$r_i$ for $p_{x_i}$ is replaced by $C_i^{-1}r_i\to0$, which is a
contradiction.  Thus $\{p_{x_i}\}$ is  normal family, and we may pass
to a subsequence such that $p_{x_i}\to p$.

Passing to further subsequences, we also have
$g_{x_i}=\psi_{x_i}(p_{x_i}(\zeta^k))\to h(p(\zeta^k))$.  Thus
$G=h(p(\zeta^k))$ has order $k^2$ at $\zeta=0$.   Since $k$ is the
maximal order on
$\cJ_k$, we have $k^2\le k$, so $k=1$.
\qed

\section 7.  Regularity

In the sequel we consider points $x$ where (\dag) holds for both $f$ and
$f^{-1}$.  We use the superscripts $u$ and $s$ to distinguish between the
``unstable''
objects $V^u$, $\Psi^u$, $E^u$, $\gamma^u$, $\tau^u$, and the ``stable''
objects
$V^s$,
$\Psi^s$, $E^s$, $\gamma^s$, $\tau^s$ (i.e.\ the corresponding objects for
$f^{-1}$).
With this notation, the backward contraction condition is now written
$\hat\lambda^u(x,n)\to0$ as $n\to-\infty$, and forward expansion is written
$\hat\lambda^u(x,n)\to\infty$ as $n\to+\infty$.  By {\it forward
contraction} we will
mean
$\hat\lambda^s(x,n)\to0$ as $n\to+\infty$, and by {\it backward expansion}
we will mean
$\hat\lambda^s(x,n)\to\infty$ as $n\to-\infty$.
$\cJ_k$ will now be written
$\cJ_{*,k}$, and the set $\cJ_j$ corresponding to $f^{-1}$ will be written
$\cJ_{j,*}$.  We set $\cJ_{j,k}=\cJ_{j,*}\cap\cJ_{*,k}$.

\proclaim Proposition 7.1.  If $V^s_x$ and $V^u_x$ exist at $x$, and if
$x$ is a point of forward expansion, then $V^s_x\ne V^u_x$, i.e.\ the germs
cannot
coincide.

\give Proof.  As was noted after (3.1), $V^u_x\subset J^-$, and
$V^s_x\subset J^+$.
Thus if $V^s_x=v^u_x$, then $V^u_x\subset\{G^+=0\}$.  If $\psi\in\Psi^u_x$ is
nonconstant, then $m_\psi(r)=\max_{|\zeta|\le r}G^+(\psi(\zeta))$ vanishes
for some
$r>0$.  But by Proposition 1.5 we cannot have $\lambda(x,n)\to\infty$ as
$n\to\infty$.  \qed

We define $\cJ_{j,k}'=\{x\in\cJ_{j,k}:\alpha(x)\cap\cJ_{j,k}\ne\emptyset\}$.

\proclaim Proposition 7.2.  If $x\in\cJ_{j,k}'$ is a point of backward
contraction,
then
$E^u_x\ne E^s_x$.

\give Proof.  Let us suppose $E^s_x=E^u_x$, and let us write $g=f^{-1}$.
Then by (4. )
$$\Vert Df^{-n}_x|_{E^s_x}\Vert^{\#,s} = \Vert
Dg^n_x|_{E^s_x}\Vert^{\#,s}\ge1$$
for $n\ge0$, where $\Vert\ \Vert^{\#,s}$ denotes the metric
$|\cdot|/\gamma^s_j$ for
$f$ (or the expanding metric for $g$), and
$$\Vert Df^{-n}_x|_{E^u_x}\Vert^{\#,u}=\lambda(x,-n).$$
Let us select a subsequence $-n_j\to-\infty$ such that $f^{-n_j}x\to\hat
x\in\cJ_{j,k}$.  By (6.2) we know that on the compact set $\{\hat
x\}\cup\{f^{-n_j}x:j=1,2,3,\dots\}$ the metrics $\Vert\ \Vert^{\#,s}$ and
$\Vert\
\Vert^{\#,u}$ are comparable to the euclidean metric, which contradicts the
backward
contraction.
\qed

We define $\tau^\iota(x)$ to be the order of intersection (or contact) between
$V^u_x$ and $V^s_x$.  Specifically, $\tau^\iota(x)=1$ means that $V^u_x$ and
$V^s_x$ meet transversally at
$x$ in the sense that $E^u_x\ne E^s_x$.  More generally, we interpret
$\tau^\iota(x)=\tau$ to mean that there is a holomorphic coordinate system
$(z,w)$ in a neighborhood of $x=(0,0)$ with
$$[V^u_{loc}(x)\cup V^s_{loc}(x)]\cap\{|z|,|w|<1\}\subset \{|z|,|w|<1,|z|\le
|w|^\tau\};\eqno(7.1)$$
and $\tau$ is the minimal value for which
this holds.  Observe that if $V^{s/u}_x$ are singular, then $\tau$ need not
be an integer.  We define
$\cJ^i_{j,k}=\{x\in\cJ_{j,k}:\tau^\iota(x)=i\}$.

If $f$ and $f^{-1}$ are quasi-expanding, then by Theorem 6.7 $\cJ_{1,1}$ is
a dense,
open subset of $J^*$. Since $\cJ_{1,*}\ni x\mapsto E^s_x$ and $\cJ_{*,k}\ni
x\mapsto E^u_x$ are continuous, and since $\cJ_{1,1}^1$ contains the saddle
points (where
$E^u_x\ne E^s_x$) it follows that $\cJ_{1,1}^1$ is a dense, open subset of
$J^*$.

It is useful to have the following quantitative version of Lemma 7.2.
\proclaim Theorem 7.3.  Let $x\in\cJ_{j,k}'$ be a point of backward
contraction.  Let
$\tau=[\tau^\iota(x)]$ denote the greatest integer in $\tau^\iota(x)$.  Then
$\tau^s(\hat x)\ge j\tau$ for all $\hat x\in\alpha(x)$.

\give Proof.  Let $(x,y)$ be a coordinate system satisfying (7.1).  For
variables
$(A,B)$ in this coordinate system we have
$$f^n(A,B)=\alpha_0(n) + \alpha_1(n)A+\dots+\alpha_{\tau-1}(n)A^{\tau-1}
+O(|A|^\tau+|B|)$$
with $\alpha_j(n)\in\C^2$.
Our first object is to show that $\alpha_r(n)\to0$ for $1\le r\le\tau-1$ as
$n\to-\infty$.  Let us choose $\psi^u\in\Psi^u_x$ such that
$$\psi^u(\zeta)=(c\zeta^k,0)+O(\zeta^{k+1})$$
for some $c\ne 0$.  If we set $\lambda=\lambda^u(\psi^u,n)$ and
$\alpha_j=\alpha_j(n)$,
then we have
$$\tilde f^n\psi^u(\zeta)=f^n\circ \psi^u(\lambda^{-1}\zeta) =
\alpha_0 + \sum_{r=1}^{\tau-1}\alpha_r\left(c
\lambda^{-k}\zeta^k+O(\lambda^{-k-1}\zeta^{k+1})\right)^r+\cdots$$
The power series coefficients of $\tilde f^n\psi^u(\zeta)$ are bounded as
$n\to-\infty$.  The coefficient of $\zeta^k$ is
$\alpha_1c\lambda^{-k}$.  By the backward contraction we have
$\lambda\to0$ as $n\to-\infty$; and since $c\ne0$, it follows
that $\alpha_1\to0$ as $n\to-\infty$.

To proceed by induction, let us suppose that $\alpha_t(n)\to0$ for $1\le t\le
r-1$ as $n\to-\infty$.  The coefficient of $\zeta^{rk}$ is
$$(\alpha_1E_1+\cdots+\alpha_{r-1}E_{r-1}+
\alpha_rc^r)\lambda^{-kr}.$$
Here the $E_t$ denote expressions in the coefficients of $\psi^u$ which are
independent of $n$.  Since $\alpha_t\to0$ as $n\to-\infty$ for $1\le t\le
r-1$, it follows that $\alpha_rc^r\lambda^{-kr}$ is bounded, so
$\alpha_r\to0$ as $n\to-\infty$.

Now we write $\psi^s(\zeta)=(\psi^s_1(\zeta),\psi^s_2(\zeta)) =
(c\zeta^j,0)+O(\zeta^{j+1})$ for some nonzero constant $c$.  If we set
$\lambda=\lambda^s(\psi^s,n)$ and $\alpha_r=\alpha_r(n)$, we have
$$\eqalign{\tilde f^n\psi^s &=f^n\circ\psi^s(\lambda^{-1}\zeta)\cr
&=\alpha_0 +\sum_{r=1}^{\tau-1}\alpha_r\left(c(\lambda^{-1}\zeta)^j
+O((\lambda^{-1}\zeta)^{j+1})\right)^r + O(|\psi^s_1|^\tau +|\psi^s_2|).\cr}$$

For $\hat x\in\alpha(x)$ there exists a sequence $n_i\to-\infty$ such that
$f^{n_i}(x)\to\hat x$.  We may pass to a subsequence so that $\tilde
f^{n_i}\psi^s$ converges to an element $\hat\psi\in\Psi^s_{\hat x}$.  Now we
have $\alpha_r\to0$ for $n\to-\infty$ for $1\le r\le\tau-1$, and (always)
$|\lambda^s_n|^{-1}\le1$, so it follows that all the coefficients of the terms
$$\alpha_r\left(c\lambda^{-j}\zeta^j +O(\lambda^{-j-1}\zeta^{j+1})\right)^r$$
tend to
zero as $n=n_i\to-\infty$.  We conclude that the only nonvanishing terms in
$\hat\psi$ arise from the expression $O(|\psi^s_1|^\tau +|\psi^s_2|)$.
However, $(\psi^s_1)^\tau=O(\zeta^{j\tau})$ by definition, and $|\psi^u_2|\le
C|\psi^u_1|^\tau = O(\zeta^{j\tau})$ by the tangency condition.  It follows
that ${\rm Ord}(\hat\psi)\ge j\tau$.
\qed
If
$V$ is a germ of a variety at $x\in\C^2$ which is locally irreducible at
$x$, then
there is a holomorphic coordinate system
$(w,z)$ in a neighborhood of $x$ such that $x=(0,0)$, and $V$ is
represented near
$(0,0)$ in terms of a Puiseaux series
$$w=a_jz^{j/m}+a_{j+1}z^{(j+1)/m}+\cdots=\sum_{n=j}^\infty a_n
z^{n/m}.\eqno(7.2)$$
(See [Ch, \S10] for details.)  Choosing the $z$-axis to be the tangent
cone, we have
$j/m>1$.  If $j/m\in\Z$, we may replace $w$ by $w'=w-a_jz^{j/m}$.  If $V$
is regular at
$x$, i.e.\ if $V$ is a complex manifold in a neighborhood of $x$, then we
may continue
this procedure and obtain a coordinate system $(w',z')$ such that
$V=\{w'=0\}$ in a
neighborhood of the origin.  If $V$ is not regular, we may continue this
procedure to
the point where we have  $a_j\ne0$, $j/m\notin\Z$, and $j/m>1$.

\proclaim Theorem 7.4.  Let $x\in\cJ'_{jk}$ be a point of backward
contraction.  Then
$V^u_x$ is a (nonsingular) manifold in a neighborhood of $x$.

\give Proof. Suppose that $V_x^u$ is not regular at $x$.  Choose a holomorphic
coordinate system $(z,w)$ at $x=(0,0)$ such that $V_x$ has a Puiseux
representation (7.1) with $a_j\ne0$, $j/m\notin\Z$, and $j/m>1$.  Now
choose $\psi\in\Psi^{k,u}_x$.
We may assume that $\psi$ has the form
$$\psi(\zeta)=(\psi_1(\zeta),\psi_2(\zeta))=
(\zeta^k+c_{k+1}\zeta^{k+1}+\cdots,\zeta^\ell+\cdots)$$
with $\ell/k=j/m$.  Let us define coefficients $\alpha_{r,t}(n)\in\C^2$
such that
$$f^n(A,B)=\sum_{r,t=1}^\infty \alpha_{r,t}(n)A^rB^t.$$
If we set $\lambda=\lambda^u(x,n)$ and $\alpha_{r,t}=\alpha_{r,t}(n)$, then
$$\tilde f^n\psi(\zeta)=f^n(\psi(\lambda^{-1}\zeta)) =
\alpha_{0,0}+\alpha_{1,0}(\lambda^{-k}\zeta^k+\cdots)
+\cdots+\alpha_{1,0}(\lambda^{-\ell}\zeta^\ell+\cdots).$$
Since $\{\tilde f^n\psi:n\le0\}$ is a normal family, all of its power series
coefficients are bounded.  The coefficient of $\zeta^k$ is
$\alpha_{1,0}\lambda^{-k}$.
Now by backward contraction it follows that $\alpha_{1,0}\to0$ as
$n\to-\infty$.

Define $q$ by the property $kq<\ell<k(q+1)$.  We next show that
$\alpha_{t,0}\to0$ as
$n\to-\infty$, for $1\le t\le q$.  We proceed by induction, assuming that
$\alpha_{t,0}\to0$ for $1\le t\le r-1$.  The coefficient of $\zeta^{rk}$ in
$\tilde
f^n\psi$ is
$$\lambda^{-rk}[\alpha_{1,0}E_{r,1} +
\alpha_{2,0}E_{r,2}+\cdots+\alpha_{r-1,0}E_{r,r-1} + \alpha_{r,0}]$$
where $E_{r,t}$ denotes a polynomial in the coefficients of $\psi$.  Now
$E_{r,t}$ is
independent of $n$, and $\alpha_{1,0},\dots,\alpha_{r-1,0}\to0$ as
$n\to-\infty$, so we
conclude that $\alpha_{r,0}\to0$.
Similarly, the coefficient of $\zeta^\ell$ is
$$\lambda^{-\ell}[\alpha_{0,1}+\alpha_{1,0}E_{\ell,1}+\cdots+\alpha_{q
,0}E_{\ell
,q}],$$
and we conclude that $\alpha_{0,1}\to0$ as $n\to-\infty$.

We conclude, therefore, that $Df^n\to0$ as $n\to-\infty$.  But let $\hat
x\in\alpha(x)\cap\cJ_{j,k}$ be given, and extract a subsequence
$-n_j\to-\infty$ such
that $f^{-n_j}\to\hat x$.  We let $g=f^{-1}$, and apply Proposition 6.1 to
$g$.  We
conclude that $\Vert D g^{n_j}|E_{g^{n_j}x}\Vert^e$ is bounded below by a
constant
(since we always have $\lambda(x,n_j)\ge1$).  Thus $Df^n$ cannot tend to
zero, and thus
$V_x$ cannot have a singular Puiseux representation.
\qed

\proclaim Corollary 7.5.  Let $x\in\cJ_k'$ be a point of backward
contraction.  Then,
modulo rotation, there is exactly one element of
$\Psi_x^{k,u}$, and this is a homogeneous parametrization.  Further, for
$\psi\in\Psi_x^{k,u}$,
$$\hat\lambda_x^k=|\lambda_\psi|^k=||Df_x||^\#.$$

\give Proof.  Let $\psi\in\Psi_x^{k,u}$ be given.  By Theorem 7.4, $V_x$ is
a regular
variety at
$x$.  Thus we may define a branch of
$\phi(\zeta):=\psi(\zeta^{1/k})$, which is holomorphic at $\zeta=0$.  Thus
$\phi'(0)\ne0$.  By Lemma 6.4,
$\psi'(\zeta)\ne0$ for
$\zeta\ne0$, so it follows that
$\phi$ is an immersion.  The uniqueness of $\psi$ now follows from Lemma
6.5.  The
equation now follows from (6.1).
\qed

We have noted earlier that $\cJ_{*,k}\ni x\mapsto E^u_x$ is continuous.  It
follows that the points of
$\cJ'_{*,k}$ are points of continuity of the metric $\cJ_{*,k}\ni
x\mapsto\Vert\
\Vert^\#_x$ on $E^u$.

\proclaim Corollary 7.6. Let $x\in\cJ'_{j,k}$ be a point of backward
contraction.
Then
$$\Vert\ \Vert^\#_x=\lim_{\cJ_{j,k}\ni y\to x}\Vert\ \Vert^\#_y.$$

\give Proof.  For $y\in\cJ_{j,k}$, we choose $\psi_y\in\Psi^{k,u}_y$ such that
$\Vert v\Vert^\#_y=|v/(\psi_y^{(k)}(0)/k!)|$.  If $\{y_i\}\subset\cJ_{j,k}$
is any
sequence converging to $x$, then $\psi_{y_i}$ converges to an element of
$\Psi_x^{k,u}$.  Since $\Psi_x^{k,u}$ consists of homogeneous paramterizations,
which are essentially unique, $\lim_{i\to\infty}\psi_{y_i}$ exists (modulo
rotation),
and thus the norms must converge.
\qed

\section 8. Hyperbolicity

In this section we explore
conditions that imply that $f$ is (uniformly) hyperbolic, as well as ways
in which
hyperbolicity can fail.  For instance in Theorem 8.3 we show that purely
geometric conditions on $J^\pm$ are sufficient to guarantee hyperbolicity.
We show in Theorem 8.4 if $f$ is
quasi-expanding, quasi-contracting and expansive, then $f$ is uniformly
hyperbolic.) Finally, we show
(Corollary 8.10) that for a special class of non-hyperbolic maps there are
points of
tangency, i.e.\ points where $E^s_x=E^u_x$.

In this section let us make the standing assumption that, unless mentioned
otherwise, $f$ is quasi-expanding and quasi-contracting.

\proclaim Proposition 8.1.  If $S\subset\cJ^i_{j,k}$ is a compact,
invariant set, then $i=1$, and $S$ is a (uniformly) hyperbolic set for $f$.

\give Proof.  Recall that $\cJ_{jk}\ni x\mapsto E^{s/u}_x$ is continuous.
By the
compactness of $S$ we have $\alpha(x)\subset S\subset\cJ^i_{jk}$ for all $x\in
S$.  Thus by Lemma 7.2 $E^s_x\ne E^u_x$, and so $i=1$.  This gives us a
continuous
splitting of $T_x\C^2$ for $x\in S$, and so by compactness the angle between
$E^s_x$ and $E^u_x$ is bounded below.  The uniform expansion/contraction of
$Df$
on $E^{s/u}$ follows from Proposition 6.1. \qed

Consider the (finite) collection of index pairs $(j,k)$ for which
$\cJ_{jk}\ne\emptyset$.  We define a partial ordering on this collection of
index pairs as follows.  We say
$(j,k)\ge(a,b)$ if $j\ge a$, $k\ge b$, and $\cJ_{jk}\ne\emptyset$.  By the
semicontinuity of $\tau^s$ and $\tau^u$,  $\cJ_{jk}$ is compact for a
maximal pair $(j,k)$.  By Proposition 8.1, then, $\cJ_{jk}$ is a
hyperbolic set for all maximal pairs $(j,k)$.

Let us consider ways in which hyperbolicity can fail to hold for $f$.  If $f$
and
$f^{-1}$ are both quasi-expanding, then  hyperbolicity (or the failure of
hyperbolicity) along an orbit is determined by the position of the orbit with
respect to the strata
$\cJ^i_{j,k}$.   For a point $x\in\cJ^1_{j,k}$, there is always uniform
contraction in the direction $E^s_x$ along the forward orbit (apply
Proposition 6.1 to $f^{-1}$).  If $Df^n|E^u_x$ is not uniformly expanding
for $n\ge0$,
then there is a subsequence $n_l\to\infty$ for which
$f^{n_l}x\to\{\tau^u>k\}$.
An alternative is that hyperbolicity may fail along a forward orbit
because the angle between $E^s_{f^nx}$ and $E^u_{f^nx}$ is not bounded
below.  In this case we have a subsequence $\{n_l\}$ with
$f^{n_l}x\to\{\tau^\iota>1\}$.  By similar reasoning, we see that the failure
of hyperbolicity along a backward orbit is caused either by
$f^{-n_m}x\to\{\tau^s>j\}$ or by $f^{-n_m}x\to\{\tau^\iota>1\}$, or both.

Theorem 8.3 gives a criterion for hyperbolicity for general polynomial
diffeomorphisms (that is to say we make no a priori assumption that $f$
is quasi-expanding and quasi-contracting) which refers only to the geometry of
$J^+$ and
$J^-$ and makes no
direct reference to $f$.  For this we will need a preliminary result.

\proclaim Lemma 8.2.  Let $\cN\subset\C^2$ be an open set, and let $\cL$ be a
Riemann surface lamination of $\cN\cap\partial K^+$.  If $T$ is a smooth
2-dimensional transversal to $\cL$ at $p\in\cN\cap\partial K^+$, then $p$ is
in the closure of $T-K^+$.

\give Proof.  Since $\cL$ is a lamination, there is a neighborhood $\cU$ of
$p$ such that the restriction $\cL|\cU$ is homeomorphic to the trivial
lamination of $S\times\Delta$, where $\Delta$ is the unit disk, and
$S\subset\Delta$ is closed.  By Slodkowski [S], $\cL$ may be extended to a
lamination $\cL^*$ of $\cU$.  Shrinking $\cU$, we may suppose that the
restriction $\cL^*|\cU$ is homeomorphic to the (trivial) lamination of
$\Delta\times\Delta$, whose leaves are $\{q\}\times\Delta$.

Since $\cL^*$ extends a lamination of $\cU\cap\partial K^+$, there are sets
$S_0,S_1\subset\Delta$ such that the leaves corresponding to
$\{q\}\times\Delta$ fill out $\cU\cap K^+$ as $q$ ranges over $S_0$ and they
fill out $\cU-K^+$ as $q$ ranges over $S_1$.  Further, the leaves
corresponding to $\partial S_0$ fill out $\cU\cap\partial K^+$.

For $x\in\Delta$, we define $\chi(x)$ as the intersection point of $T$ and the
leaf corresponding to $\{x\}\times\Delta$.  Since $T$ is transversal, $\chi$
is defined and continuous (possibly after shrinking $\cU$).  Let $\hat
p\in\partial S_0$ be such that $\chi(\hat p)=p$.  Now there are points $q\in
S_1$ arbitrarily close to $\hat p$, and so the points $\chi(q)\in T-K^+$ are
arbitrarily close to $p$.\qed

In the following theorem we make no a priori assumption about quasi-expansion
or quasi-contraction.

\proclaim Theorem 8.3.  A polynomial diffeomorphism of $\C^2$ is hyperbolic
on $J^*$
  (resp.\ $J$) if and only if there is a neighborhood $\cN$ of $J^*$ (resp.\
$J$),
and Riemann surface laminations $\cL^\pm$ of $\cN\cap J^\pm$ such that $\cL^+$
and $\cL^-$ intersect transversely at all points of $J^*$ (resp.\ $J$).

\give Proof. We start by working with $J^*$.  The fact that this lamination
structure exists for a hyperbolic set of a diffeomorphism is standard. We will
prove the converse. For a saddle point
$p$, it follows from  (7) of [BS6, Theorem 2.1] that $W^u(p)$ is a leaf of
$\cL^-$.
The lamination hypothesis implies that the leaves of $\cL^-$ may be written
locally
as a family of graphs of holomorphic functions.  Since bounded analytic
functions
have locally bounded first derivatives, this implies that the bounded area
condition
(3.1) holds.  For each $p\in J^*$, whether or not $p$ is a saddle, the variety
$V^u(p,\epsilon)$ is a manifold which is transversal to $\cL^+$ at $p$.  By
Lemma 8.2, $V^u(p,\epsilon)$ intersects $\C^2-K^+$ arbitrarily close to
$p$. In particular
the function $G^+$ is positive on $V^u(p,\epsilon)$. Compactness of the set of
varieties
$V^u(p,\epsilon)$ gives a positive lower bound for the maximum of
$G^+$ on $V^u(p,\epsilon)$ which is independent of $p$. Thus by Theorem 3.4,
$f$ is quasi-expanding. By Proposition 5.3, $\cJ_1=J^*$. By similar arguments,
$f^{-1}$ is quasi-expanding, and $J^*=\cJ_{1,1}$.  By Proposition 8.1,
then, $J^*$
is a hyperbolic set for $f$.

Now let us deal with $J$.  Under our hypotheses the currents $\mu^\pm$
supported
on $J^\pm$ are given by transverse measures.  Thus the wedge product, $\mu$,
can be interpreted locally as a product measure. It follows that the support of
$\mu$, which is a priori a subset of $J$, is actually equal to $J$. But the
support of $\mu$ is $J^*$. According to the previous paragraph, $J=J^*$ is a
hyperbolic set for $f$. \qed

Let $\psi=(\psi_1,\psi_2)=(z(\zeta),w(\zeta)):(\C,0)\to(\C^2,x)$ be a germ of a
holomorphic mapping, and let $V(\psi)$ denote the induced germ at $x$.
Then $V(\psi)$
has a Puiseux representation (7.2) with $j>m$, so that $V(\psi)$ is tangent
to the
$z$-axis.  Since $V(\psi)$ is locally irreducible at $x$, we may assume
that $\gcd(m,n_1,n_2,\dots)=1$, where $n_i$ is a listing of all the numbers
such that
$a_{n_i}\ne0$.

Let us recall some facts about complex varieties.  (See [Ch, \S10,\S12] for
details).
For a point $y\in V$, we let $\mu(V,y)$ denote the multiplicity of
$V$ at $y$.  This number is defined by the property that for a generic
complex line $L$
passing near $y$, $L\cap V$ contains exactly $\mu(V,y)$ points near $y$.
$V$ is regular
at $y$ if and only if $\mu(V,y)=1$.  If $V$ is written as a Puiseux
expansion (7.2),
with $\gcd(m,n_1,n_2,\dots)=1$, then any line $L$ transversal to $\{w=0\}$
and  passing near $(0,0)$ will intersect $V$ in $m$ points near the origin.
Thus
$\mu(V,x)=m$.

If ${\rm Ord}(\psi)=k$, then ${\rm Ord}(\psi_1)=k\ge m$, and ${\rm
Ord}(\psi_2)=kj/m>k$.  Let $\xi(\zeta)=\zeta+\cdots$ be a germ of a holomorphic
function such that $\psi_1(\zeta)=c\xi^k$ near $\zeta=0$.  We may assume
that $c=1$.
We may write $\psi_2(\zeta)=\sum c_n\xi^n$.  This gives another Puiseux
representation
for $V$: $w=\sum c_n z^{n/k}$.  On the other hand, the Puiseux
representation is
essentially unique.  So $m$ divides $k$, and we may set $p:=k/m\in\Z$.
Thus we may
write $\psi(\zeta)=\phi(\xi^p)=(\phi_1(\xi^p),\phi_2(\xi^p))$, where
$\phi_1(t)=t^{k/p}$, and $\phi_2(t)=\sum_{n=j}^\infty a_nt^{n/p}$, where $n$ is
divisible by $p$ whenever $a_n\ne0$.

To summarize, if ${\rm
Ord}(\psi)=k>m$, then $m$ divides $k$, and $\psi$ covers the variety $V$
exactly $p=k/m$
times.  The relation between the multiplicity
(order) of the parametrizing function and the multiplicity of the variety
is thus
$${\rm Ord}(\psi)=p\cdot\mu(V,x).\eqno(8.1)$$

In the sequel we will treat $V(\psi)$ as the variety $V$, but counted with
multiplicity $p$.  One reason for introducing multiplicities is that it
makes it
easier to view varieties as currents: if
$\phi_j$ is a sequence of nonsingular mappings which converge to
$\phi$ in some neighborhood of the origin, then the corresponding germs
$V(\phi_j)$
converge as currents (in some neighborhood of $x$) to the current defined by
$V(\psi)$ counted with multiplicity $p$.

If $V_1$ and $V_2$ are 1-dimensional varieties which intersect only at $x$,
we may
define $\iota_x(V_1,V_2)$, the intersection number at $x$.  This number has the
property that for almost every small $\tau_1,\tau_2\in\C^2$ the translates
$V_j+\tau_j$, $j=1,2$, intersect in $\iota_x(V_1,V_2)$ points near $x$.  In
general,
we have
$$\iota_x(V_1,V_2)\ge\mu_x(V_1)\mu_x(V_2).$$
Equality holds if the tangent cones of $V_1$ and $V_2$ at $x$ are distinct.
The
intersection number behaves continuously: if
$V_1^j$  (resp.\ $V_2^j$) are sequences of varieties that converge in the
sense of
currents to
$V_1$ (resp.\
$V_2$), then for $j\ge j_0$, we have
$$\iota_x(V^j_1,V^j_2)=\iota_x(V_1,V_2).$$

Now let
$\psi_j$, $j=1,2$ be a germs of a mapping that define the varieties $V_j$,
$j=1,2$.  If
$x$ is an isolated point of intersection of $V_1$ and $V_2$, then
$$ \eqalign{\iota_x(V(\psi_1),V(\psi_2))
&=p_1p_2\cdot\iota_x(V_1,V_2)\cr
&\ge p_1p_2\cdot\mu_x(V_1)\mu_x(V_2) = {\rm Ord}(\psi_1){\rm
Ord}(\psi_2)}
\eqno(8.2)$$
These properties of varieties give us the following:

\proclaim Lemma 8.4.  Let $\psi_j:\{|\zeta|<1\}\to\C^2$, $j=1,2$ be nonconstant
mappings with $\psi_1(0)=\psi_2(0)$.  Set $m_j={\rm Ord}(\psi_j)$, and let
$r>0$ be
given.  Then for $\hat\psi_j$ sufficiently close to $\psi_j$, there are sets
$X_j\subset\{|\zeta|<r\}$ such that
$$\sum_{a\in X_j}\iota(V_1,V_2,\hat\psi_j(a))\cdot{\rm Ord}(\hat\psi_j,a)
\ge m_j,$$
where ${\rm Ord}(\hat\psi_j,a)={\rm Ord}(\hat\psi_j(\zeta-a))$ denotes the
order of
$\hat\psi_j$ at $\zeta=a$.

An important topological dynamical consequence of hyperbolicity is
the shadowing property. The following result gives us
a quantitative measure of the failure of uniqueness of shadowing.

\proclaim Theorem 8.5.  Suppose that $f$ is quasi-expanding and
quasi-contracting.  If
$x\in\cJ_{j,k}$, then for $\epsilon>0$ there is a set $X\in J^*$ containing
$jk$ elements such that
$$\sup_{n\in\Z}\max_{a,b\in X}dist(f^na,f^nb)<\epsilon.$$

\give Proof.  Choose $\psi^s\in\Psi^s_x$ with ${\rm Ord}(\psi^s)=j$ and
$\psi^u\in\Psi^u_x$ with ${\rm Ord}(\psi^u)=k$.  If $\{p_i\},\{q_j\}\subset
S$ are
sequences converging to $x$ with $\phi_{p_i}\to\psi^s$ and
$\phi^u_{q_i}\to\psi^u$,
then by Lemma 6.4 of [BLS], we may assume that $\phi^s_{p_i}(|\zeta|<1)$
(resp.\ $\phi^s_{p_i}(|\zeta|<1)$) intersects
$\psi^u(|\zeta|<1)$ (resp.\ $\psi^s(|\zeta|<1)$) transversally.

Let $M=\sup_{\psi\in\Psi}\max_{|\zeta|\le1}|\psi'(\zeta)|$, and set
$r=\epsilon/M$. Let
$X^s,X^u\subset\{|\zeta|<r\}$ be the sets given by Lemma 8.4.  Since
$\phi^s_{p_i}$ and
$\phi^u_{q_i}$ are immersions, the order at each point is equal to 1.  And
since the
immersions are transversal for $i$ sufficiently large, the intersection
numbers are 1.
Thus each set $X^s$ and $X^u$ contains at least $jk$ points.

Let $X=\{\phi^s_{p_i}(\zeta):\zeta\in X^s\}=\{\phi^u_{q_i}(\zeta):\zeta\in
X^u\}$.
Thus for $n\ge0$ we have
$$\max_{a,b\in X}dist(f^na,f^nb) =\max_{\zeta',\zeta''\in X^s}dist
(f^n\phi^s_{p_i}(\zeta'),f^n\phi^s_{p_i}(\zeta''))$$
$$=\max_{\zeta',\zeta''\in X^s}
dist(\phi^s_{p_i}(\lambda'_n\zeta'),\phi^s_{p_i}(\lambda_n''\zeta''))
\le\max_{|\zeta|\le\kappa^{-n}r}|\phi^s_{p_i}(\zeta)'|\le\kappa^{-n}rM
\le\epsilon.$$
For $n\le0$, we use $\phi^u_{q_i}$ instead, and we conclude that the
diameter of $f^nX$
is no greater than $\epsilon$ for all $n\in\Z$.  \qed

Note that  if the varieties $V^s_x$ and $V^u_x$ are tangent, then the set $X$
may be taken to have strictly more than $jk$ elements.

\proclaim Corollary 8.6.  Suppose that $f$ is quasi-expanding and
quasi-contracting,
but $f$ is not hyperbolic on $J^*$.  Then $f$ is not expansive.

\give Proof.  If $f$ is not hyperbolic, it follows from Proposition 8.1 that
$\cJ_{j,k}\ne\emptyset$ for some index pair $(j,k)\ne(1,1)$.  By Theorem
8.5, then, $f$
is not expansive.  \qed
Since hyperbolic mappings are expansive, and expansivity is preserved under
topological
conjugacy, we have the following.
\proclaim Corollary 8.7.  If $f$ is quasi-expanding and quasi-contracting
but not
hyperbolic, then $f$ is not topologically conjugate to a hyperbolic map.

We will use the following result, which is a special case of Proposition
5.1 of [V].
\proclaim Lemma 8.8.  Suppose that $V$ is a subvariety of the bidisk
$\{|z|,|w|<1\}$,
and suppose that the projection to the $z$-axis is proper and has degree
bounded by
$m<\infty$.  For any $\epsilon>0$ there is a $\delta>0$, depending only on
$\epsilon$
and $m$, such that if
$\tilde V$ is a connected component of $V\cap\{|z|<\delta\}$, then the
diameter of
$\tilde V$ is less than $\epsilon$.

We refer to $\cC:=J^*-\cJ_{1,1}$ as the singular locus of $f$.  In the
following
results, we consider
$f$ for which $\cC$ is finite.  This is parallel to the
critical finiteness condition in one complex dimension.  Note that if $\cC$ is
finite, then
$\cC$ consists of saddle points, and $\cV^{s/u}$ are regular on $J^*$ and form
laminations on $\cJ_{1,1}=J^*-\cC$.  Further, $E^s$ and $E^u$ are transverse at
$\cC$.  Thus the set of tangencies, written $\cT=\{x\in J^*:
E^s_x=E^u_x\}$, is a
subset of $\cJ_{1,1}=J^*-\cC$.

\proclaim Proposition 8.9.  Let $f$ be quasi-expanding and
quasi-contracting, and let
$\cC$ be finite.  Then for each tangency $r\in\cT$, there are points
$p,q\in\cC$ such
that $r\in W^s(p)\cap W^u(q)$.

\give Proof.  If $r\in\cT$, then $\alpha(r)\subset\cC$ by Theorem 7.3. Thus
$r\in
W^u(q)$ for some $q\in\cC$.  Similarly, $\omega(r)\subset\cC$, so $r\in W^s(p)$
for some $p\in\cC$.\qed

\proclaim Theorem 8.10.  If $f$ is quasi-expanding and quasi-contracting,
and if
$\cC$ is finite, then $\cT$ is a discrete subset of $\cJ_{1,1}$, and the
closure
of $\cT$ is $\cT\cup\cC$.

\give Proof.  Since $\cC=J^*-\cJ_{1,1}$ is finite, it consists of periodic
points,
which must be saddle points by Corollary 6.3.  Saddle points are
not points of tangency, so $\cT\subset\cJ_{1,1}$.
The families of varieties $\cV^s$ and
$\cV^u$ are laminations in a neighborhood of $\cJ_{1,1}$.  Thus any
tangency must be
isolated by Lemma 6.4 of [BLS].

Now let us fix a saddle point $x\in \cC$; and passing to a higher iterate
of $f$,
we may assume it is a fixed point.  It follows that $x\in\cJ_{j,k}$ for some
index pair $(j,k)\ne(1,1)$.  We may assume that $k>1$.  Consider a coordinate
system $(z,w)$ such that $x=(0,0)$, and $f$ is essentially linear on
$\cB:=\{|z|,|w|<1\}$, with uniform expansion in the horizontal direction and
uniform contraction in the vertical direction.  Given $\rho>0$, we may choose
small $r>0$ and $0<\rho_1<\rho_2<1$ such that for any point $q$ of
$\{|z|<1,|w|<r\}$ there is an $n>0$ such that $f^{-n}q$ belongs to
$S:=\{|z|<\rho,\rho_1<|w|<\rho_2\}$.  For
$\rho,\rho_2$ and
$\cB$, let $\cV^s(\cB)$ and $\cV^u(\cB)$ denote the set of varieties in $\cB$
corresponding to $V^{s/u}(x)$ for $x\in\{|z|<\rho,|w|<\rho_2\}$.  We will
choose
$\rho,\rho_2$ and $\cB$ small enough that for
$x\in J^*\cap\{|z|<\rho,|w|<\rho_2\}$, $V^u_x$ is a subvariety of
$\{|z|<1,|w|<1\}$, with proper projection to the $z$-axis and a uniform
bound on
the mapping degree of the projection.

Let us choose $\psi\in\Psi^{k,u}_x$, and let $\{p_i\}\subset S$
be a sequence such that $\phi_{p_i}$ converges to $\psi$.  For $i$ sufficiently
large, we have $V(p_i)\subset\{|z|<1,|w|<r\}$, and we may choose the first
$n_i$ such that $f^{-n}V^u(p_i)\cap S\ne\emptyset$.  Let us choose a
subsequence of
$\{\tilde f^{-n_i}\psi_{p_i}\}$ which converges to a limit
$\tilde\psi\in\Psi_y$ for
some $y$ in the closure of $S$.  Since $V^u(y)\ne V^u(x)$, it follows that
$y\in\cJ_{*,1}$.

Given $\epsilon>0$, let us choose $\delta>0$ as in Lemma 8.8.  For $i$
sufficiently
large, $\tilde V_i:= V^u(p_i)\cap\{|z|<\delta\}$ is connected.  Since $n_i$
was chosen
for the first time $f^{-n_i}V_i$ intersects $S$, it follows that
$f^{-n_i}\tilde V_i\subset\cB$.  Since $f$ is contracting in the
$z$-direction, it
follows that the projection of
$f^{-n_i}\tilde V_i$ to the $z$-axis has diameter less than $\delta$.  By
Lemma 8.8,
then the diameter of $f^{-n_i}\tilde V_i$ is less than $\epsilon$.

To see that $V(y)$ intersects $W^s_{loc}(p)$ tangentially, we recall that
$V^u(\psi_{p_i})\cap V^s(x)\cap\cB$ consists of $k$
points, which are also contained in $f^{-n_i}\tilde V_i$.  By Lemma 8.8,
the set of
intersection has diameter no greater than $\epsilon>0$.  Since these points
remain
inside a compact subset of
$\cB$, it follows that the intersection multiplicity of $V^u(y)$ and
$W_{loc}^s(p)$ is
also
$k$, and the diameter of the set of intersection is $\le\epsilon$.  Since
$\epsilon$
may be taken arbitrarily small, $V^u(y)$ intersects $W^s_{loc}(p)$ in a
single point
of multiplicity $k$.  Since $k>1$, this is a tangency.
\qed
Taking into account the multiplicity $k$ in the last paragraph of this proof,
we have the following.
\proclaim Corollary 8.11.  Suppose that $f$ is quasi-expanding and
quasi-contracting
but is not hyperbolic.  If $\cC$ is finite, then $\cC\ne\emptyset$, and
there are
points of tangency.  More precisely, if $p\in\cJ_{j,k}$, $k>1$, then there
is a point
$y\in W^s(p)$ where
$W^u(y)$ is tangent to $W^s(p)$, and the order of contact is $k$.
Conversely, if
$y\in W^s(p)$ is a point of tangency between $W^s(p)$ and $W^u(y)$, then
the order of
contact is no greater than
$k$.

\section A. Appendix:  One-dimensional Mappings

In this paper we have developed an approach to the study of a dynamically well
behaved family of maps of $\C^2$ via a family of immersions from $\C$ into
$\C^2$.  In this Appendix we explore a similar approach to one-dimensional
mappings.  Our purpose is to justify the analogy between semi-hyperbolicity in
$\C$ and quasi-hyperbolicity in $\C^2$.  This is achieved in Theorem A.5.

Let us consider a  polynomial mapping $g:\C\to\C$ of
degree
$d>1$.  Let $J=J_g$ denote the Julia set, and let $K=K_g$ denote the filled
Julia set, so
$J=\partial K$.  Let $G=G_K$ denote the Green function of $\C-K$ with pole
at infinity.
Let $S\subset J$ denote the set of repelling periodic points.  For each
$x\in S$ we let $n$ denote the period of $x$, so that $g^nx=x$.  We define
$\lambda(x,n):=g^n(x)'$ and
$L_n(\zeta)=\lambda(x,n)\zeta$.  There is a (linearizing) function
$\phi_x:\C\to\C$ such
that
$$\phi_x(0)=x, {\rm\ \ and\ \ }g^n\circ\phi_x(\zeta)
=\phi_x(\lambda(x,n)\zeta) = \phi_x\circ L_n(\zeta)\eqno(A.1)$$
(see Milnor [M]).  The linearizing function is the analog of the unstable
manifold, and the functional equation is the analogue of (1.4).  The function
$\phi_x$ also satisfies
$$\phi_x=g^n\circ\phi_x\circ L_n^{-1}=g^{jn}\circ\phi_x\circ
L^{-j}_n\eqno(A.2)$$
for every $j\ge0$.   If $\phi_x'(0)=1$, we may define $\phi_x$ simply as
$$\phi_x(\zeta)=\lim_{j\to\infty} g^{nj}\circ L^{-j}_n.\eqno(A.3)$$
For $k\ge0$, $g^k\phi_x:(\C,0)\to(\C,g^kx)$ is a linearizing function at
$g^kx$.
For $\alpha\in\C$, $\alpha\ne0$,  $\phi_x(\alpha\zeta)$ is also a linearizing
function.  We fix $0<t<\infty$, and we define $\psi_x$ to be the
linearizing function
$\psi_x:\zeta\mapsto\phi_x(\alpha\zeta)$, with $|\alpha|$ determined by
condition
(1.13).  Thus we have a family of maps
$\psi_S=\{\psi_x:\C\to\C:x\in S\}$.

As in \S1 we may take
normal limits and obtain the family $\Psi$, where each $\psi\in\Psi_x$ is
defined and
holomorphic on a domain $\Omega_x$ with
$\{|\zeta|<1\}\subset\Omega_x\subset\C$.  We may
define the transformation  $\tilde g:\Psi_x\to\Psi_{gx}$ as in \S1, and if
$\psi_x$ is
nonconstant, we may define the multiplier $\lambda=\lambda_{\psi_x}$ by the
condition
$\tilde g(\psi_x)(\zeta)=\psi_{gx}(\lambda^{-1}\zeta)$.

We will say that $g$ is {\it quasi-expanding} if $\Psi$ is a normal family of
entire functions.   By Proposition 1.7, quasi-hyperbolicity is independent of
normalizing constant $t$; it will be convenient for us to choose a specific
value
of $t$ just before Lemma A.4. By Proposition 1.5, quasi-expansion implies that
(\ddag) holds at each $x\in J$.  By Theorem 1.2, it is equivalent to
$|\lambda_x|\ge\kappa>1$ for all $x\in S$.

If $g$ is quasi-expanding, we define $\tau$ as was done just before Proposition
5.1.  There is a natural stratification
$J=\cJ_1\sqcup\cdots\sqcup \cJ_k$, where $\cJ_m=\{\tau=m\}$.  We define the
infinitesimal metric $\Vert\cdot\Vert^\#_x$ on the tangent space $T_x\C$
for $x\in
\cJ_i$ as in \S6.  Note that $x\mapsto\|\cdot\|_x^\#$ is not globally
continuous,
but it is continuous on each stratum $\cJ_i$.  This metric is uniformly
expanded
by $g'$.  And as in Corollary 6.3, there is a $\kappa>1$ such that
$|g^n(x)'|>\kappa^n$ holds for each point of period $n$.  While it was known
earlier that a semi-hyperbolic map has a (singular) metric which is uniformly
expanded (see [Ca]), this construction for quasi-expanding maps seems more
direct, in addition to defining an (infinitesimal) metric at each point of $J$.

Let $\cC=\{z\in\C:g'(z)=0\}$ denote the set of critical points of $g$.  For
$c\in\cC$, let $\cP(c)=\{g^j(c):j\ge1\}$, let
$\cP(\cC)=\bigcup_{c\in\cC}\cP(c)$,
and let $\bar\cP(c)$ denote the closure of $\cP(c)$.

\proclaim Lemma A.1.  If $c\in\cC$ be a critical point, then $\tau\ge2$ on
$\bar\cP(c)$.
If $c\in\omega(c)$, then $\tau=\infty$ on $\omega(c)$.

\give Proof.  If $\psi\in\Psi_c$ is constant, then $\tilde
g^j(\psi)\in\Psi_{g^jc}$
is constant.  Thus $\tau(g^jc)=\infty$ for all $j\ge0$.  Since $\tau$ is upper
semicontinuous, it is equal to $\infty$ on the closure of
$\{g^j(c):j\ge1\}$.  Now let
$\psi\in\Psi_c$ be a nonconstant function.  It follows that $\tilde
g^j(\psi)$ has a
critical point at the origin for $j\ge1$, i.e.\ $\tau>1$ on
$\{g^j(c):j\ge1\}$.  Again,
by upper semicontinuity, $\tau>1$ on the closure of this set.

Now suppose that $c\in\omega(c)$. If $\Psi_c$ consists only of the constant
function,
then $\tau(c)=\infty$.  If there is a nonconstant $\psi\in\Psi_c$, then ${\rm
Ord}(\psi)<\infty$.  Let
$n_j\to\infty$ be a sequence such that
$g^{n_j}c\to c$.  By the chain rule, ${\rm Ord}(\tilde g^{n_j}\psi)>{\rm
Ord}(\psi)$.
Passing to a subsequence of $\{n_j\}$, we may assume that $\tilde
g^{n_j}\psi\to\hat\psi$.  By the upper semicontinuity of $\tau$, we have
${\rm Ord}(\hat\psi)>{\rm Ord}(\psi)$.  Thus
$\tau(c)=\infty$.   \qed

We let $\cC'_x:=\{\zeta\in\C:\psi_x'(\zeta)=0\}$ denote the set of critical
points of
$\psi_x$.
\proclaim Lemma A.2.  If $x\in S$ is a repelling periodic point, then
$\psi_x(\cC'_x)\subset\cP(\cC)$.

\give Proof. Suppose $\zeta\in\cC_x'$.  Since $x$ is a repelling periodic
point,
$|\lambda(x,n)|>1$, so it follows from (A.1) that
$\psi_x'(0)\ne0$.  Let $U$ be a neighborhood of the origin in $\C$ where
$\psi_x'\ne0$.  Choose $j$ such that $\tilde\zeta= L^{-j}\zeta\in U$, and set
$\tilde z=\psi_x(\tilde\zeta)$.  By (A.2) and the Chain Rule,
$$\psi_x'(\zeta)=(g^{jn}\circ\psi_x\circ L^{-j}(\zeta))' =
g'(g^{jn-1}(\tilde z))\cdots g'(g(\tilde z))\cdot g'(\tilde
z)\cdot(\psi_x\circ L^{-j})'(\zeta)=0.$$ It follows $g'(g^k\tilde z)=0$ for
some $1\le
k\le jn-1$, which means that
$g^k\tilde z\in\cC$.  Thus $x=g^{jn-k}(g^k\tilde z)\in\cP(\cC)$.
\qed
\proclaim Proposition A.3.  If $g$ is quasi-expanding, then
$\cJ_1=J-\bar\cP(\cC)$.

\give Proof.  By Lemma A.1, $\cJ_1$ is disjoint from $\bar\cP(\cC)$.
Conversely, let
$y\notin\bar\cP(\cC)$ be given.  Choose $0<\delta<dist(y,\bar\cP(\cC))$.  Let
$x$ be a repelling periodic point sufficiently close to $y$ that
$B(x,\delta)\cap\bar\cP(\cC)=\emptyset$.  By Lemma A.2, there are no
critical values of
$\psi_x$ in the disk $B(x,\delta)$.  Thus there is an analytic function
$\phi:B(x,\delta)\to\C$ such that $\psi\circ\phi(z)=z$.  By the Koebe
Distortion
Theorem, $\{|\zeta|<\delta|\phi'(x)|/4\}\subset \phi(B(x,\delta))$.

Now let
$\chi(x)=\inf\{|\zeta|\in\cC'_x\}$.  Since
$\phi(B(x,\delta))\cap\cC_x'=\emptyset$,
we have $\delta|\phi'(x)|/4\le\chi(x)$, or
$\delta/4\le\chi(x)|\psi_x'(0)|$.  Since $|\psi'_x(0)|$ is bounded above,
it follows that
$\chi(x)$ is bounded below.

To show that $y\in J_1$, we need to show that $\psi_y'(0)\ne0$ for every
$\psi_y\in\Psi_y$.  Let us take a sequence  $x\to y$, such that
$\psi_x\to\psi_y$.  Since
$\chi(x)$ is bounded below, there is an open neighborhood $U$ of the origin
in $\C$ where $\psi_x'\ne0$ on $U$ for all $x$.
The limit $\psi_y'$ is then either nonvanishing on $U$, or it vanishes
identically.
By (\ddag), then,  $\psi_y$ does not vanish on $U$. \qed

For a domain $D$ and $y\in g^{-n}D$, we let $(g^{-n}D)_y$ denote the
connected component
of
$g^{-n}D$ containing $y$.  A mapping
$g$ is  said to be {\it semi-hyperbolic} (see [CJY]) if there are numbers
$\epsilon_0>0$ and
$M<\infty$ such that for every
$n\ge0$,
$0<\epsilon<\epsilon_0$ and $x\in J$, the mapping degree of
$$g|(g^{-n}B(x,\epsilon))_y:(g^{-n}B(x,\epsilon))_y\to B(x,\epsilon)$$
is bounded by $M$
for each $y\in g^{-n}x$.  If $g$ is semi-hyperbolic, then by Theorem
2.5 there is an
$a>0$ such that for $\rho_1>0$ sufficiently small, depending only on $M$,
such that for
all $x\in J$, all $n\ge0$, and all $y\in g^{-n}x$, we have
$$B(y,as)\subset (g^{-n}B(x,\rho_1\epsilon))_y\subset B(y,s)$$
for some $s>0$.  Set
$$t:=\min_{x\in J}\max_{B(x,\rho_1\epsilon)}G.$$
\proclaim Lemma A.4.  Let $g$ be semi-hyperbolic, and let
$\epsilon,\rho_1,t,a>0$ and
$M<\infty$ be  as above.  Then there is a number $B<\infty$ such that for
any periodic
point $x$,
$\omega_x:=(\psi_x^{-1}B(x,\rho_1\epsilon))_0$ satisfies
$$\{|\zeta|<B^{-1}\}\subset\omega_x\subset\{|\zeta|<B\}$$
and $\psi_x:\omega_x\to B(x,\epsilon)$ is a proper mapping of degree $\le M$.

\give Proof.  Let $n$ denote the period of $x$, and assume that $x=0$.  For
$0<\rho<1$,
let $\bar t:=\max_{x\in J}\max_{B(x,\rho\rho_1\epsilon)}G$.  Choose $\rho$
small enough
that $\bar t<t$.

For each $j\ge0$ we have $\{|\zeta|<ar\}\subset\omega\subset \{|\zeta|<r\}$ for
some $r=r_j$ corresponding to
$\omega=(L^j_ng^{-jn}B(0,\rho\rho_1\epsilon))_0$ and for
some $r=\tilde r_j$ corresponding to
$\omega=(L^j_ng^{-jn}B(0,\rho_1\epsilon))_0$.
We may take the limit as $j\to\infty$ in (A.3) so that $g^{jn}\circ
L^{-j}_n\to\phi_x$,
and we may pass to a subsequence to have $r_j\to r$ and $\tilde r_j\to
\tilde r$.  If we
write $\tilde\omega_x'=(\phi^{-1}_xB(0,\rho_1\epsilon))_0$ and
$\omega_x'=(\phi^{-1}_xB(0,\rho\rho_1\epsilon))_0$, then we have
$\{|\zeta|<r\}\subset\omega_x'\subset\{|\zeta|<r\}$ and
$\{|\zeta|<\tilde r\}\subset\tilde\omega_x'\subset\{|\zeta|<\tilde r\}$.  Thus
$$\tilde\omega_x'-\bar\omega_x'\subset\{ar<|\zeta|<\tilde r\}.$$
Thus the moduli satisfy
$$\log(\tilde r/(ar))\le {\rm Mod}(\{ar<|\zeta|<\tilde r\})\le{\rm
Mod}(\tilde\omega'_x-\bar\omega'_x).$$

Let us remark that $\phi_x:(\phi^{-1}_xB(0,\epsilon))_0\to B(0,\epsilon)$
is a proper
mapping with degree bounded by $M$ since each $g^{jn}\circ L^{-j}_n$ was
also a proper
mapping with degree bounded by $M$.  It follows that
$\phi_x:(\tilde\omega'_x-\bar\omega'_x)\to (B(x,\rho_1\epsilon)-\bar
B(x,\rho\rho_1\epsilon))$ is a proper map.  The modulus of an annulus is
defined as
the extremal length of the family of curves connecting the two boundaries (cf.\
[A, Chapter 4]).  Under a proper map, this family pulls back to a family of
curves
which connect the two boundaries; thus the modulus cannot decrease, so we have
$${\rm Mod}(\tilde\omega'_x-\bar\omega'_x)\le {\rm
Mod}(B(x,\rho_1\epsilon)-\bar
B(x,\rho\rho_1\epsilon))=\log(1/\rho).$$
We conclude that $\tilde r/(ar)\le\rho^{-1}$.

Finally, let us consider $\psi_x$ and $\omega_x$, which are obtained from
$\phi_x$ and
$\tilde\omega_x'$ by a scaling by a linear factor $\lambda>0$.  Thus
$\{|\zeta|<\lambda
ar\}\subset\omega_x=\lambda\tilde\omega'_x\subset\{|\zeta|<\lambda\tilde
r\}$.  By the
definition of $t$, we have $\max_{B(x,\rho_1\epsilon)}G=\max_{\omega_x}G\ge
t$.  By the
Maximum Principle, we have
$\max_{|\zeta|<\lambda\tilde r}G\ge t$.  It follows by (1.13) that
$\lambda\tilde r\ge1$.
Similarly, we have
$\max_{B(x,\rho\rho_1\epsilon)}G=\max_{\lambda\omega'_x}G\le\bar t<t$.
Again by the Maximum Principle, $\max_{|\zeta|<\lambda ar}G\le t$.
Thus by (1.13) we have  $ar\lambda<1$.  By our previous inequality, it
follows that
$\rho\le\lambda ar\le \lambda\tilde r\le\rho^{-1}$, so we may take
$B=\rho^{-1}$.

Finally, the mapping degree of the restriction of $\psi_x$ to
$(\psi_x^{-1}B(0,\epsilon))_0$ is the same as the degree of restriction of
$\phi_x$, so
it is bounded by $M$. \qed
We will use the following estimate on the Green function  (see [CJY, \S3]):
{\sl If  $g$
is semi-hyperbolic, then there exist $\eta>0$ and $A<\infty$ such that
$$\max_{B(x,r)}G\ge \eta r^A \eqno(A.4)$$
for all $x\in J$, $0<r<1$.}  Note that (1.12) and (A.4) are similar but
different; the
estimate (A.4) takes place on dynamical space while (1.12) concerns the
uniformizations.

\proclaim Theorem A.5.  Quasi-expansion $\Leftrightarrow$ semi-hyperbolicity.

\give Proof.  Suppose first that $g$ is quasi-expanding.  By Corollary 6.3, any
periodic point is expanding.  Thus there are no parabolic points.  Now
suppose that
$c\in J$ is a critical point.  By quasi-expansion, we must have
$\tau<\infty$ on $J$,
so by Lemma A.1 this means that $c$ is not contained in $\omega(c)$, its
$\omega$-limit set.  It follows by [CJY, Theorem 1.1] that $g$ is
semi-hyperbolic

Now suppose that $g$ is semi-hyperbolic.  Let $B$ be as in Lemma A.4, and
choose $\chi>B^2$.  By Theorem
2.5, we may choose
$\rho_2>0$ sufficiently small that for any $y\in J$ there is a number
$s=s_y$ such
that
$$(\psi_y^{-1}(B(y,\rho_2\rho_1\epsilon))_0\subset
\{|\zeta|<s\}\subset\{|\zeta|<\chi
s\}\subset(\psi_y^{-1}(B(y,\rho_1\epsilon))_y.$$
By the right-hand inclusion in Lemma A.4, we have $\chi s<B$.

Let us set
$$\bar t_1:=\max_{x\in J}\max_{B(x,\rho_1\epsilon)}G.$$
Choose $k$ such that
$$\left({\bar t_1 d^{-k}\over\eta}\right)^{1/A}\le a\rho_2\rho_1\epsilon.$$
By Theorem 2.5, there is an $r>0$ such that
$B(y,ar)\subset (g^{-k}(B(x,\rho_1\epsilon))_y\subset B(y,r)$
for $y\in g^{-k}x$.  By the definition of $\hat t_1$ and the maximum
principle, we have that $G\le \hat t_1d^{-k}$ on $B(y,ar)$.  It follows
from (A.4) that
$$\bar t_1d^{-k}\ge \eta(ar)^A.$$
By the choice of $k$ we conclude that $r\le \rho_2\rho_1\epsilon$.
Thus
$$\psi_y^{-1}(g^{-k}(B(x,\rho_1\epsilon))_0
\subset(\psi_y^{-1}B(x,\rho_2\rho_1\epsilon))_0\subset\{|\zeta|<s\}$$
with $s$ as above.

Let
$L:\zeta\mapsto\lambda\zeta$ denote the linear map such that
$g^k\circ\psi_y=\psi_x\circ L$.  By this functional equation, $L$ maps
$(\psi_y^{-1}g^{-k}B(x,\rho_1\epsilon))_0\subset\{|\zeta|<s\}$ to
$(\psi_x^{-1}(B(x,\rho_1\epsilon))_0$.  This last set contains
$\{|\zeta|<1/B\}$ by the
left-hand containment in Lemma A.4.  Thus
$\lambda\ge (sB)^{-1}$, which is no smaller than $\chi B^{-2}$ since $\chi
s<B$.
We conclude that $|\lambda|$ is uniformly bounded below by $\kappa:=\chi
B^{-2}>1$, so by
Theorem 1.2, $g^k$ is quasi-expanding.  By Proposition 1.3, then, $g$ is
quasi-expanding.
\qed
A consequence of Corollary 6.3 is:
\proclaim Corollary A.6.  If $g$ is semi-hyperbolic, then the repelling
periodic
points are uniformly repelling.

Questions dealt with in this Appendix also arise naturally in connection
with the
study of the structure of leaves in the induced inverse limit system.  See [LM]
for this approach.
\bigskip

\centerline{\bf References}

\item{[A]}  L. Ahlfors, {\sl Conformal Invariants}, McGraw-Hill, New York,
1973.

\item{[BLS]}  E. Bedford, M. Lyubich, and J. Smillie, Polynomial
diffeomorphisms of $\C^2$.  IV: The measure of maximal entropy and laminar
currents,  Invent.\ Math.\ 112, 77--125 (1993).

\item{[BS1]}  E. Bedford and J. Smillie,  Polynomial diffeomorphisms of
$\C^2$: Currents, equilibrium measure, and hyperbolicity, Invent. Math. 87,
69--99 (1990).

\item{[BS6]}  E. Bedford and J. Smillie,  Polynomial diffeomorphisms of
$\C^2$.  VI: Connectivity of $J$,  Annals of Math.\ 148, 695--735 (1998).

\item{[BS7]}  E. Bedford and J. Smillie,  Polynomial diffeomorphisms of
$\C^2$.  VII: Hyperbolicity and external rays,  Ann.\ Sci.\ Ecole Norm.\ Sup.,
32, 455--497 (1999).

\item{[BS]} E. Bedford and J. Smillie, Real polynomial diffeomorphisms with
maximal entropy: Tangencies.

\item{[CJY]}  L. Carleson, P. Jones, and J.-C. Yoccoz, Julia and John,
Bol.\ Soc.\ Bras.\ Mat., Vol.\ 25, 1--30 (1994).

\item{[Ca]}  J. Carette, Liens entre la g\'eometrie et la dynamique des
ensembles
de Julia, Th\`ese, U. de Paris-Sud, Orsay, 1997.

\item{[Ch]} E. Chirka, {\sl Analytic Sets},  Kluwer Academic Publishers, 1989.

\item{[F]}  W.H.J. Fuchs, {\sl Topics in the Theory of Functions of One Complex
Variable}, van Nostrand, 1967.

\item{[HO]} J.H.\ Hubbard and R.\ Oberste-Vorth, H\'enon mappings in the
complex domain II: Projective and inductive limits of polynomials, in: {\sl
Real and Complex Dynamical Systems}, B. Branner and P. Hjorth, eds. 89--132
(1995).

\item{[LS]} F. Ledrappier and J.-M. Strelcyn, A proof of the estimation from
below in Pesin's entropy formula, Ergodic Theory Dyn.\ Syst.\ 2, 203--219
(1982).

\item{[LM]} M. Lyubich and Y. Minsky, Laminations in holomorphic dynamics,
J. of
Differential Geometry, 47, 17--94 (1997).

\item{[M]} J. Milnor, Dynamics in One Complex Variable: Introductory
Lecture, Stony
Brook IMS Preprint 1990\#5.

\item{[S]} Z. Slodkowski, Holomorphic motions and polynomial hulls, Proc.\
Amer.\ Math.\ Soc.\ 111, 347--355 (1991).

\item{[V]}  K. Verma, Boundary regularity of correspondences in $\C^2$,
Math.\ Z., 231,
253--299 (1999).


\rightline{Eric Bedford}
\rightline{Indiana University}
\rightline{Bloomington, IN 47405}
\bigskip
\rightline{John Smillie}
\rightline{Cornell University}
\rightline{Ithaca, NY 14853}

\bye